\documentclass[12pt]{article}
\usepackage{amssymb,amsmath,amsthm,epsfig}
\usepackage{cite}
\usepackage{eepic}
\usepackage{epic}
\usepackage{graphicx}
\usepackage{color}
\usepackage{ifpdf}
\usepackage[title]{appendix}
\usepackage{amsthm}

\usepackage{dsfont}
\usepackage{graphicx} 
\usepackage{subfigure}
\usepackage{algorithm}
\usepackage{algorithmic}

\usepackage{amsmath}
\usepackage{bm}
\usepackage{multirow}

\usepackage{caption}

\usepackage{tikz}
\usetikzlibrary{shapes.geometric,arrows}
\usepackage{makecell}

\usepackage{graphicx}
\usepackage{epstopdf}
\usepackage[colorlinks=true,
linkcolor=blue,
anchorcolor=blue,
citecolor=blue,
urlcolor=blue,]{hyperref}
\usepackage{cleveref}
\usepackage{notoccite}

\theoremstyle{plain}
\newtheorem{lem}{Lemma}[section]
\newtheorem{thm}[lem]{Theorem}

\theoremstyle{definition}

\theoremstyle{remark}
\newtheorem{rem}{Remark}[section]

\begin{document}
	
	\title{\large\bf {{Model free data assimilation with Takens  embedding}}}
	\author{Ziyi Wang\thanks{School of Mathematical Sciences,  Tongji University, Shanghai 200092, China. ({\tt  2111156@tongji.edu.cn}).}
		\and
		Lijian Jiang\thanks{School of Mathematical Sciences,  Tongji University, Shanghai 200092, China. ({\tt  ljjiang@tongji.edu.cn}).}
	}
	\date{}
	\maketitle
	\begin{abstract}
	    In many practical scenarios, the dynamical system is not available and standard data assimilation methods are not applicable. Our objective is to construct a data-driven model for state  estimation without the underlying dynamics. Instead of directly modeling the observation operator with noisy observation, we establish  the state space model of the denoised observation. Through data assimilation techniques, the denoised observation information could be used  to recover the original model state. Takens' theorem shows that an embedding of the partial and denoised observation is diffeomorphic to the attractor. This gives   a theoretical base  for estimating the model state using the reconstruction map. To realize the idea, the procedure consists of  offline stage and online stage.  In the offline stage,  we construct the surrogate dynamics using dynamic mode decomposition with noisy snapshots to learn  the transition operator for the denoised observation. The filtering distribution of the denoised observation can be estimated using adaptive ensemble Kalman filter, without knowledge of the model error and observation noise covariances. Then the reconstruction map can be established using the posterior mean of the embedding and its corresponding state. In the online stage, the observation is filtered  with the surrogate dynamics. Then the online state estimation can be performed utilizing the reconstruction map and the filtered observation. Furthermore, the  idea can be generalized to the nonparametric framework with nonparametric time series prediction methods for chaotic problems. The numerical results show the proposed method can estimate the state  distribution without the physical  dynamical system.
	\end{abstract}\smallskip
	
	\smallskip
	
	{\bf keywords:}
	data assimilation, dynamic mode decomposition, Takens' theorem

	\section{Introduction}
	Data assimilation (DA) aims to estimate model state for a dynamical system, combining partial and noisy observation with prior information via Bayes’s formula \cite{D-A}. DA has broad   applications in various real-world scenarios, such as global positioning systems\cite{GPS}, oil recovery\cite{oil-recory}, and weather forecasting\cite{weather-forecasting}. Developing practical  data assimilation methods for state estimation is of practical significance and has attracted great attention from scientists and engineers.
	
	The most popular data assimilation method is  Kalman filter\cite{kalman-filter}, proposed by Rudolf Kalman, which provides an exact algorithm to determine the filtering distribution for linear problems with additive Gaussian noise. The basic idea of Kalman filter can be generalized to address nonlinear and non-Gaussian problems. These methods such as extended Kalman filter\cite{ExKF} and emsemble Kalman filter\cite{EnKF}, update the model state by  a Gaussian assumption. The extended Kalman filter realizes state estimation  by linearizing the dynamics. The ensemble Kalman filter (EnKF) generates an ensemble of particles and updates the ensemble by the linear assumption. As approximating Gaussian filters, they may not offer an accurate approximation of the true filtering distribution; nevertheless, they exhibit robustness employing in high dimension\cite{D-A}. In contrast, as a sequential Monte Carlo method, the particle filter\cite{pf} is able to estimate the true filtering distribution in the limit of a large number of particles. However, the particle filter  does  not perform well in practice because particle degeneracy usually occurs for high-dimension problems. All these classical filtering methods require a specific form of the underlying dynamics, and their performance is significantly  contingent upon the underlying dynamical system  and model assumptions. In many practical applications, the explicit dynamical system is often  unavailable. This  motivates us to combine data-driven modeling methods with data assimilation\cite{kalman-net,em-da,modeling-1,modeling-2} for the cases without any physical dynamics.
	
	Data-driven modeling methods have emerged as promising approaches applied to DA, such as dynamic mode decomposition (DMD)\cite{DMD-1,DMD-2,DMD-3}, sparse identification of nonlinear dynamical systems (SINDy)\cite{SINDy}, etc. When  the transition model is unknown, these techniques offer an effective surrogate for the transition operator and enable  the successful execution of DA within the underlying dynamic system \cite{DMDEnKF}. When the observation operator is identical, the method proposed by\cite{EnKF-DMD} can effectively denoise the model state by offering a proxy transition model with DMD. Kalman-Takens method \cite{takens-filter} models the transition operator in the scenarios where only partial state data is available. This approach is based on  Takens' theorem \cite{takens, takens-nonscalar}, which shows that a time-delayed embedding of partial and denoised observation is diffeomorphic to the attractor of a compact manifold. Furthermore, the research  in \cite{takens-similar} has shown that the mapping from a time-delayed embedding to the attractor can be effectively learned to reconstruct the true model state using supervised learning techniques, such as neural networks (NNs). However, certain limitations exist within the current research. For instance, the Kalman-Takens method \cite{takens-filter} is unable to estimate the unobservable state variables. The methodology proposed in \cite{takens-similar} requires noise-free observation data, which limits its applicability in the scenarios where observation data is inherently noisy.
	
	In this paper, we propose a  method, DMD-Takens based EnKF (DMD-T EnKF), which  facilitates real-time estimation with data-driven modeling methods in the case where the underlying dynamical system and noise information are unknown.  Considering the challenges in directly modeling the observation operator, we establish the state space model of the denoised observation. The denoised observation is then estimated with DA. This allows  us to estimate the model state utilizing the embedding of the denoised observations through a reconstruction map based on  Takens' theorem. In the offline stage, we firstly construct the surrogate transition model for the denoised observation by DMD with noisy observations. Given the transition model, the filtering distribution of the denoised observation is estimated by adaptive ensemble Kalman filter (adaptive EnKF)\cite{adaptive-EnKF}, without knowledge of the model error and observation noise covariances. The posterior mean of the denoised observation will be leveraged to construct a more accurate transition model. Besides, the filtering distribution of the time-delayed embedding can be approximated with the denoised observations. Consequently, we learn the reconstruction map through supervised learning methods, utilizing the posterior mean of the embedding and its corresponding model state data. In the online stage, we estimate the denoised observation by adaptive EnKF and construct  time-delayed coordinate vector. After that, the filtering distribution of model state is approximated with posterior samples obtained from the reconstruction map. Furthermore, our method can be extended  to a nonparametric framework: KNN-Takens based EnKF (KNN-T EnKF). The transition model in the denoised observation space, along with the reconstruction map is learned using nonparametric time series prediction methods for chaotic problems\cite{nonparametric,nonparametric-2,analog-EnKF}. Compared with Kalman-Takens method\cite{takens-filter}, our method can estimate all the model state through the reconstruction map.
	
	The paper is organized as follows. Section \ref{section2} provides a brief overview of EnKF. In Section \ref{section3}, we present the construction of surrogate model in the offline stage and provide an online inference approach. In Section \ref{section4}, a nonparametric framework is introduced for online state estimation. In Section \ref{section5}, we present several numerical examples to approximate the filtering distribution without the underlying dynamics. The paper concludes with a summary in Section \ref{section6}.

	\section{Preliminaries} \label{section2}
	Let us consider the dynamical system governed by the map $\Phi \in C(\mathbb{R}^n,\mathbb{R}^n)$ with noisy observation $y=\{y_k\}_{k \in N}$ depending on the observation operator $h\in C(\mathbb{R}^n,\mathbb{R}^m)$,
	
	\begin{equation}
	\begin{cases}
		x_{k+1}= \Phi (x_k) +\xi,\quad x_0 \sim N(m_0,C_0),\\
		y_{k+1}=h(x_{k+1})+\eta,\quad k\in N,	
	\end{cases}
	\label{ssm}
	\end{equation}
	where $\xi \sim \mathcal{N}(0,Q)$ and $\eta \sim \mathcal{N}(0,R)$. Data assimilation is concerned with two fundamental issues: smoothing and filtering. In this paper, our focus is specifically on the filtering problem. We firstly introduce EnKF\cite{EnKF}, which can effectively determine the filtering distribution for nonlinear and non-Gaussian problems by a Gaussian assumption.
	
	Let $Y_k:=\{y_\tau\}_{\tau=1}^k $ represent the accumulated observation data up to time step $k$. EnKF comprises two essential steps: prediction and analysis for evaluating filtering distribution $p(x_{k-1} | Y_{k-1}) \rightarrow p(x_{k}| Y_{k})$. Let $ \{x_{k-1}^{a,1}, \cdots, x_{k-1}^{a,N}\}$ be the samples drawn from the filtering distribution $p(x_{k-1}|Y_{k-1})$, which is approximated by
	\[
	p(x_{k-1}|Y_{k-1})\approx \frac{1}{N} \sum_{n=1}^{N} \delta_{(x_{k-1}^{a,n})} (x_{k-1}),
	\]
	where $\delta$ is the Dirac delta function and $N$ denotes the ensemble size. For the prediction step, we can apply the transition model to obtain the ensembles and approximate the prediction distribution at time $k$:
	\begin{equation}
		\begin{aligned}
			p(x_k|Y_{k-1})&=\int p(x_k,x_{k-1}|Y_{k-1})d x_{k-1}\\
			&=\int p(x_k|x_{k-1}) p(x_{k-1}|Y_{k-1})d x_{k-1}\\
			&\approx \frac{1}{N} \sum_{n=1}^{N} \mathcal{N}(x_k|\Phi (x_{k-1}^{a,n}),Q)\\
			&\approx \frac{1}{N} \sum_{n=1}^{N} \delta_{(x_{k}^{f,n})} (x_{k}),
			\label{enkf-1}
		\end{aligned}
	\end{equation}
	where $x_{k}^{f,n}=\Phi(x_{k-1}^{a,n})+\xi_k$ represents the samples drawn from the prediction distribution $p(x_{k}|Y_{k-1})$.
	
	In the analysis step, according to Bayes' rule, the filtering distribution at time $k$ satisfies
	\[
	p(x_k|Y_k)\propto p(y_k|x_k)p(x_k|Y_{k-1})\approx\frac{1}{N}\sum_{n=1}^{N} \mathcal{N}(y_k|h(x_k),R) \mathcal{N}(x_k|\Phi (x_{k-1}^{a,n}),Q).
	\]
	If the observation operator is linear, the mixture of Gaussians $p(x_k|Y_{k-1})$ in $(\ref{enkf-1})$ is approximated by a single Gaussian distribution. Consequently, the posterior distribution can be approximated by a Gaussian distribution,
	\[
	p(x_k|Y_k)\propto p(y_k|x_k) p(x_k|Y_{k-1}) \approx \mathcal{N}(y_k|Hx_k,R) \mathcal{N}(x_k|x_k^f,\hat{\Sigma}_{k}),
	\]
	where $H$ is the linear observation operator,
	$x_{k}^f=\frac{1}{N}\sum_{n=1}^{N}x_k^{f,n}$, $\hat{\Sigma}_{k}=\frac{1}{N-1}\sum_{n=1}^N (x_k^{f,n}-x_{k}^f)(x_k^{f,n}-x_{k}^f)^{T}$. The main steps of EnKF read as following,
	\[
	\begin{aligned}
		&\hat{K}=\hat{\Sigma}_{k}H^{T}(H\hat{\Sigma}_{k}H^{T}+R)^{-1},\\
		&x_k^{a,n}= x_k^{f,n}+\hat{K}(y_k-y_k^{n}),\\
		&y^n_k=H x_k^{f,n}+\eta_{k}^{n},\\
		&x_{k}^a=\frac{1}{N}\sum_{n=1}^{N}x_k^{a,n},
	\end{aligned}
	\]
	where $x_k^a$ denotes the posterior mean and $y_k^{n}$ represents a pseudo-observation with $\eta_k^{n}\sim \mathcal{N}(0,R)$. If the noise statistics $Q$ and $R$ are unavailable, adaptive EnKF\cite{adaptive-EnKF} utilizes the cross-correlations of innovations to estimate them. This method can realize state estimation in the absence of noise statistics.
	
	We can find that EnKF relies on the explicit form of the transition operator $\Phi$ and observation operator $h$. In many practical scenarios, the dynamical system (\ref{ssm}) is unavailable and standard algorithms are not applicable. Instead, we can obtain trajectories from historical observations or numerical simulations. Therefore, without the underlying dynamics, our goal is to realize online state estimation using a labeled dataset $\{x_k^{(i)},y_k^{(i)}\}_{k=1}^T, i=1,2,\cdots,M$, where $k$ denotes the time index and $i$ represents the $i$-th trajectory.

	\section{DMD-Takens based EnKF}\label{section3}
	In this section, we describe DMD-Takens based EnKF, which can provide state estimation without knowledge of the underlying dynamics. Based on Takens' theorem, we estimate the denoised observation with DA and learn the reconstruction map via a supervised learning method. This method includes two stages: the offline stage and the online stage. We will present the basic framework of the proposed approach and the details of  the two stages.

	\subsection{Basic framework}
	Instead of directly modeling the original dynamics as specified in (\ref{ssm}), we establish a state space model of the denoised observation. The model state $x_k$ can be estimated through the reconstruction map derived from Takens' theorem.

	In our context, the denoised observation
	is represented by $\hat{y}_k$. Assuming that the evolution process of $\hat{y}_k$ is Markovian, the discrete-time dynamical system for $\hat{y}_k$ is defined by
	\begin{equation}
		\begin{cases}
			\hat{y}_{k+1}=f(\hat{y}_k) + \omega_k,\\
			y_{k+1}=\hat{y}_{k+1}+\eta_{k+1},
		\end{cases}
		\label{ssm_f}
	\end{equation}
	where $f$ is the transition model, and $\omega=\{\omega_k\}_{k \in N}$ represents model error, which is assumed as an i.i.d sequence. In this section, we employ DMD to construct the surrogate transition model for $\hat{y}_k$. DMD\cite{DMD-1,DMD-2,DMD-3} is a data-driven modeling method and numerically approximates the Koopman operator in a finite-dimensional space. We assume that this space is spanned by a set of observables, represented as $\Psi=\{\psi_1,\psi_2,\cdots,\psi_m,\cdots,\psi_{N_k}\}$. The first $m$ components $\{\psi_1,\cdots,\psi_m\}$ is assumed as $m$ scalar-valued functions containing the full-state observable $g(\hat{y})=\hat{y}$, i.e.,
	\[
	\psi_i(\hat{y})=\mathbf{e}_i^{T}\hat{y},\quad i=1,2,\cdots,m,
	\]
	where $\mathbf{e}_i$ is the $i$-th unit vector in $\mathbb{R}^m$. Then we define the data matrices $\hat{Y}_0$ and $\hat{Y}_1$ as follows,
	\[
	\begin{aligned}
		\hat{Y}_0 &:=[\Psi(\hat{y}_1),\cdots,\Psi(\hat{y}_{l-1})],\\
		\hat{Y}_1 &:=[\Psi(\hat{y}_2),\cdots,\Psi(\hat{y}_l)],
	\end{aligned}
	\]
	where $l$ represents the serial number of snapshots. DMD aims to compute the matrix $K$, which approximates the Koopman operator in the sense of least-squares,
	$$
	K = \arg \min _{\tilde{K}} \Vert \hat{Y}_1- \tilde{K} \hat{Y}_0\Vert_{F},
	$$
	where $\Vert \cdot \Vert_F$ represents the Frobenius norm. The best-fit matrix can be computed  by $K=\hat{Y}_1 \hat{Y}_0^{\dagger}$, where $\dagger$ is the Moore-Penrose inverse. The transition model $f$ is then approximated by the first $m$ rows of $K$, denoted as $K^{y}$. Consequently, the dynamics of $\hat{y}_k$ can be approximated as follows,
	\begin{equation}
		\begin{cases}
			\hat{y}_{k+1}=K^{y} \Psi(\hat{y}_{k})+\hat{\omega}_k,\\
			y_{k+1}=\hat{y}_{k+1}+\eta_{k+1},
		\end{cases}
		\label{dmd-ssm1}
	\end{equation}
	where $\hat{\omega}_k= f(\hat{y}_k)-K^{y} \Psi(\hat{y}_{k}) +\omega_k$ includes two primary components: the approximation error of $f(\hat{y}_k)-K^{y} \Psi(\hat{y}_{k})$ and model error $\omega_k$, which is assumed to satisfy $\hat{\omega}_k\sim \mathcal{N}(0,M
	)$. The calculation of $K^{y}$ will be described in the following subsection. If $K^{y}$ is available, the filtering distribution $p(\hat{y}_k|Y_k)$ can be estimated with adaptive EnKF, which can adaptively estimate the noise statistics $M$ and $R$. Then we aim to evaluate the filtering distribution $p(x_k|Y_k)$ of the model state through the denoised observation.
	\begin{rem}
		The transition model is constructed using DMD or  SINDy.
	\end{rem}
	
\begin{figure}[htp]
		\centering
		\includegraphics[scale=0.6]{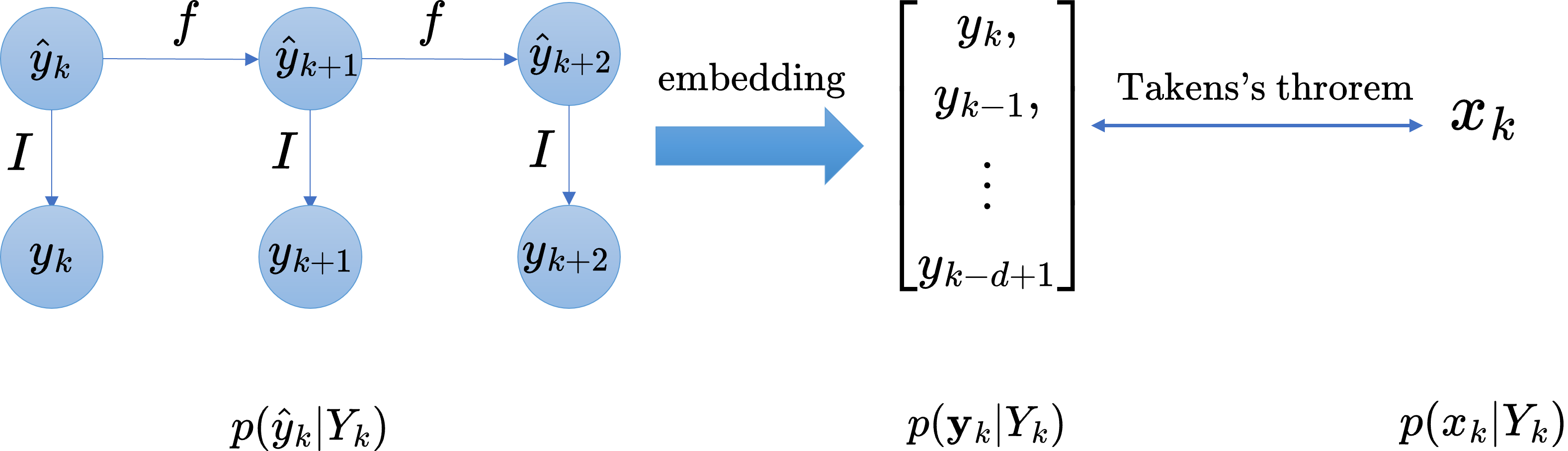}
		\caption{The general framework for DMD-Takens based EnKF.}
		\label{shili1}
	\end{figure}

	When the observation  is partial, the information in $\hat{y}_k$ is inadequate to accurately  estimate the model state $x_k$.
	\begin{thm}{[Takens, 1981]}
		Let $M$ be a compact manifold of dimension $m$. For pairs $(\varphi, g)$, where $\varphi: M \rightarrow M$ is a smooth (at least $C^2$) diffeomorphism and $g: M \rightarrow \mathbb{R}$ a smooth function, it is a generic property that the $(2 m+1)-$ delay observation map $\mathcal{G}_{(\varphi, g)}: M \rightarrow \mathbb{R}^{2 m+1}$ given by
		$$
		\mathcal{G}_{(\varphi, g)}(z)=\left(g(z), g \circ \varphi(z), \ldots, g \circ \varphi^{2 m}(z)\right)
		$$
		is an embedding.
	\end{thm}
	This theorem can be extended to the case of vector observations\cite{takens-nonscalar}. If $\varphi=\Phi^{\dagger},\ g=h$, we have
	\[
	\mathcal{G}_{(\varphi, g)}(x_k)=\left(\hat{y}_k,\hat{y}_{k-1}, \ldots,\hat{y}_{k-2m}\right).
	\]
	Based on Takens' theorem, we can construct the embedding $\mathcal{G}_{(\varphi, g)}(x_k)$ and estimate the model state $x_k$ via a reconstruction map. For convenience, we denote the embedding as $\mathbf{y}_k:=\left(\hat{y}_k,\hat{y}_{k-1}, \ldots,\hat{y}_{k-d+1}\right)$, where $d$ represents the time-delay length and can be determined by the false nearest neighbors methods\cite{FNN,FNN-nonscalar}. The reconstruction map is defined as:
	$$
	x_k=F(\mathbf{y}_k,\theta),\quad F: \mathbb{R}^{d \times m} \rightarrow  \mathbb{R}^{n},
	$$
	where $\theta$ represents the parameters of this map. We can estimate $p(x_k|Y_k)$ using samples drawn from $p(\mathbf{y}_k|Y_k)$ via the map $F$. Moreover, if $p(\mathbf{y}_k|Y_k)$ can be accurately evaluated, we  can accurately estimate the  filtering distribution $p(x_k|Y_k)$ with a large number of samples, i.e.,
	$$
	\begin{aligned}
		p(x_k|Y_k)
		&= \frac{\int p(x_k,\mathbf{y}_k,Y_k)\mathrm{d}\mathbf{y}_k }{p(Y_k)}\\
		&= \frac{\int p(x_k|\mathbf{y}_k)p(\mathbf{y}_k,Y_k)\mathrm{d}\mathbf{y}_k}{p(Y_k)}\\
		&= \int p(x_k|\mathbf{y}_k)p(\mathbf{y}_k|Y_k)\mathrm{d}\mathbf{y}_k\\
		&\approx \sum_{j=1}^{N^{'}} p(x_k|\mathbf{y}_k^{j})\\
		&= \sum_{j=1}^{N^{'}} \delta_{F(\mathbf{y}_k^j,\theta)}(x_k),
	\end{aligned}
	$$
	where $\mathbf{y}_k^{j}$ represents the $j$-th sample drawn from $p(\mathbf{y}_k|Y_k),j=1,2,\cdots,N^{'}$. We will show that the filtering distribution $p(\mathbf{y}_k|Y_k)$ can be approximated by the distribution $p(\hat{y}_k|Y_k)$ in the following section. The general framework of our approach is described in Fig.\;\ref{shili1}. We aim to model the surrogate transition model $K^{y}$ and the reconstruction map $F$ for online estimation.
	
	\subsection{Offline data-driven modeling}\label{exactDMD}
	In this stage, we learn the surrogate transition model $K^{y}$ and the reconstruction map $F$ with the labeled dataset $\{x_k^{(i)},y_k^{(i)}\}_{k=1}^T$.
	
	Inspired by EnKF-DMD in \cite{EnKF-DMD}, the observation operator can be considered as identical. The DMD matrix $K^{y}$ can be initialized using noisy snapshots $\{y_1^{(i)},\cdots,y_T^{(i)}\}_{i=1}^{M}$, which are denoted by $K^{y,(0)}$. For simplicity of notations,   the subscript $i$ will be suppressed in the rest of paper because  the method  keeps the same for  each trajectory. Accordingly, the adaptive EnKF is employed to estimate the filtering distribution $p(\hat{y}_k|Y_k,K^{y,(0)})$ and generate the posterior ensemble $\{\hat{y}_k^{a,1},\cdots,\hat{y}_k^{a,N}\}$, thereby facilitating noise reduction in the observation data. Subsequently, the posterior mean data set, denoted as $\{\hat{y}_1^{a},\cdots,\hat{y}_T^{a}\}$, is utilized to construct a more refined transition model $K^{y,(1)}$. In practice, this process can be applied iteratively to derive $K^{y,(t)}, t\geq 2$. This iterative process will be elucidated from the perspective of the expectation-maximization (EM) algorithm in Section \ref{em_section}.
	
	Then we estimate the delay coordinate vector $\mathbf{y}_k$ through the estimates of $\{\hat{y}_k,\cdots,\hat{y}_{k-d+1}\}$. The filtering distribution of $p(\mathbf{y}_k|Y_k)$ satisfies
	\begin{equation}
		\begin{aligned}
			p(\hat{y}_{k},\cdots,\hat{y}_{k-d+1}|Y_k)
			&= \underbrace{p(\hat{y}_{k-d+1}|Y_k)}_{\text{\uppercase\expandafter{\romannumeral1}}}  \prod_{i=0}^{d-2} \underbrace{ p(\hat{y}_{k-i}|\hat{y}_{k-i-1},Y_k)}_{\text{\uppercase\expandafter{\romannumeral2}}},\\
		\end{aligned}
		\label{pdf}
	\end{equation}
	where $\text{\uppercase\expandafter{\romannumeral1}}$ and $\text{\uppercase\expandafter{\romannumeral2}}$ represent the conditional distribution of $\hat{y}_{k-d+1}$ and $\hat{y}_{k-i}$ given the data sets $\{Y_{k-d+1},y_{k-d+2},\cdots,y_{k}\}$ and $\{\hat{y}_{k-i-1},Y_{k-i},y_{k-i+1}\cdots,y_{k}\}$, respectively. We assume that $Y_{k-d+1}$ and $\{\hat{y}_{k-i-1},Y_{k-i}\}$ retain sufficient information for reliable analysis. Therefore, the distributions $p(\hat{y}_{k-d+1}|Y_{k-d+1})$ and $p(\hat{y}_{k-i}|\hat{y}_{k-i-1},Y_{k-i})$ are employed to approximate $\text{\uppercase\expandafter{\romannumeral1}}$ and $\text{\uppercase\expandafter{\romannumeral2}}$, separately.
	
\begin{figure}[htp]
		\hspace{-5pt}
		\includegraphics[scale=0.8]{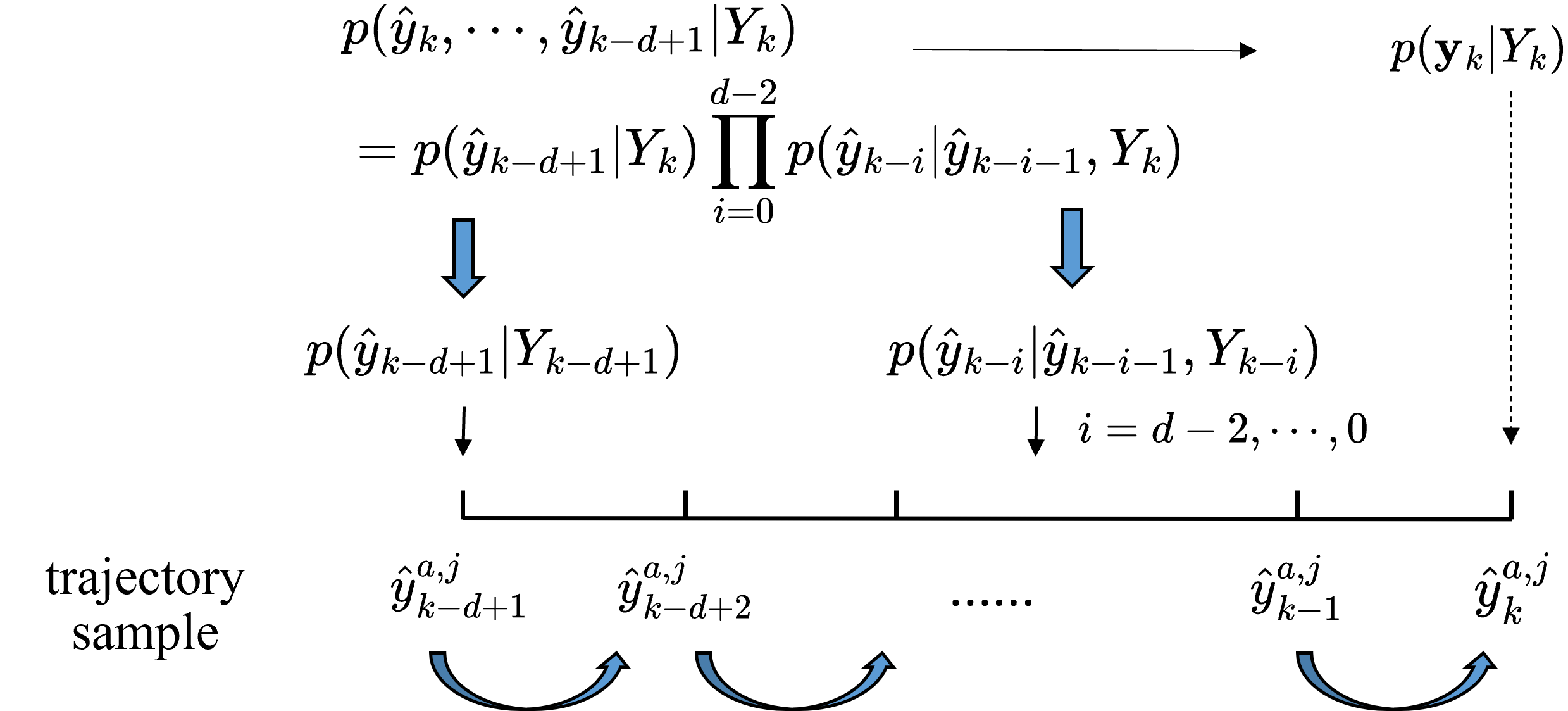}
		\caption{The sampling process for $p(\mathbf{y}_k|Y_k)$.}
		\label{shili2}
	\end{figure}

	To sample from $(\ref{pdf})$, we firstly sample from the filtering distribution $p(\hat{y}_{k-d+1}|Y_{k-d+1})$. The posterior ensemble $\{\hat{y}_{k-d+1}^{a,1},\cdots,\hat{y}_{k-d+1}^{a,N}\}$ is available from the previous process. Then we sample from $p(\hat{y}_{k-d+2}|\hat{y}_{k-d+1}^{a,j},Y_{k-d+2}), j=1,2,\cdots,N$, which satisfies
	$$
	\begin{aligned}
		p(\hat{y}_{k-d+2}|\hat{y}_{k-d+1}^{a,j},Y_{k-d+2})&=p(\hat{y}_{k-d+2}|\hat{y}_{k-d+1}^{a,j},y_{k-d+2})\\
		&\propto p(\hat{y}_{k-d+2}|\hat{y}_{k-d+1}^{a,j}) p(y_{k-d+2}|\hat{y}_{k-d+2}).
	\end{aligned}
	$$
	The particle $\hat{y}_{k-d+2}^{a,j}$ is  generated using adaptive EnKF and  can be considered as a sample  from the distribution $p(\hat{y}_{k-d+2}|\hat{y}_{k-d+1}^{a,j},Y_{k-d+2})$. We can recursively sample from the distribution  $p(\hat{y}_{k-i}|\hat{y}_{k-i-1}^{a,j},Y_{k-i}),\ \ i=d-3,\cdots,0$. The trajectory
	\begin{equation}
		\mathbf{y}_k^{a,j}:=\{\hat{y}_{k-d+1}^{a,j},\cdots,\hat{y}_{k}^{a,j}\}
		\label{trajectory}
	\end{equation}
	can be regarded as a sample drawn from the joint pdf $p(\mathbf{y}_k|Y_k)$, i.e.,
	\[
	p(\mathbf{y}_k|Y_k)\approx\frac{1}{N} \sum_{j=1}^{N} \delta_{\mathbf{y}_k^{a,j}} (\mathbf{y}_k).
	\]
	The sampling process of $p(\mathbf{y}_k|Y_k)$ is depicted in Fig.\;\ref{shili2}.

	The posterior mean $\mathbf{y}_k^{a}$ can be seen as maximum a posterior (MAP) estimation of the truth $\mathbf{y}_k$. We can utilize the dataset $\{x_k^{(i)},\mathbf{y}_k^{a,(i)}\}_{k=1}^T$ to learn the reconstruction map $F$ with supervised learning methods. This map can be approximated with parametric methods, such as fully connected neural networks, radial basis function network, etc.  The following is the loss function to learn the reconstruction map $F$,
	\[
	\mathcal{L}_F=\frac{1}{MT}\sum_{i=1}^{M}\sum_{k=1}^{T}\Vert x_k^{(i)} -F(\mathbf{y}_k^{a,(i)},\theta) \Vert _2 ^2.
	\]
	 We can also learn the reconstruction map with nonparametric methods, which will be described in the following section.
	The general framework to construct the surrogate model is described in Algorithm $\ref{algorithm1}$.
	
	\begin{algorithm}[h]
		\caption{DMD-Takens based EnKF: construction of the surrogate model}
		\textbf{Input}: labeled dataset $\{x_k^{(i)},y_k^{(i)}\}_{k=1}^T, i=1,2,\cdots,M$\\
		\textbf{Output}: transition model $K^{y}$, reconstruction map $F$\\
		~1:~ Estimate $K^{y}$ and the distribution $p(\hat{y}_k|Y_k)$. 
		The initial DMD matrix $ K^{y,(0)}$ is calculated by $\{y_1^{(i)},\cdots,y_T^{(i)}\}_{i=1}^{M}$. \\
		(i): calculate $p(\hat{y}_k^{(i)}|Y_k,K^{y,(t)})$ and obtain the posterior ensemble $\{\hat{y}_{k}^{a,1,(i)},\cdots,\hat{y}_k^{a,N,(i)}\}$ with adaptive EnKF.\\
		(ii): update the DMD matrix $K^{y,(t+1)}$. Let
		$$
		\begin{aligned}
			&Y_{0}^{a}:=[\hat{y}_{1}^{a,(1)},\cdots,\hat{y}_{T-1}^{a,(1)},\cdots,\hat{y}_{1}^{a,(M)},\cdots,\hat{y}_{T-1}^{a,(M)}],\\
			&Y_{1}^{a}:=[\hat{y}_{2}^{a,(1)},\cdots,\hat{y}_{T}^{a,(1)},\cdots,\hat{y}_{2}^{a,(M)},\cdots,\hat{y}_{T}^{a,(M)}],\\
		\end{aligned}
		$$
		we can obtain the matrix $K^{(t+1)}=\Psi(Y_{1}^{a}) \Psi(Y_{0}^{a})^{\dagger}$. Then $K^{y,(t+1)}$ is available with the first $m$ rows of the DMD matrix $K^{(t+1)}$. \\
		~2:~ Construct the sample $
		\mathbf{y}_k^{a,j,(i)}:=\{\hat{y}_{k-d+1}^{a,j,(i)},\cdots,\hat{y}_{k}^{a,j,(i)}\}
		$, i.e.,
		$$
		p(\mathbf{y}_k^{(i)}|Y_k)=\frac{1}{N} \sum_{j=1}^{N} \delta_{\mathbf{y}_k^{a,j,(i)}} (\mathbf{y}_k).
		$$
		~3:~ Learn the reconstruction map $F$ utilizing the dataset $\{x_k^{(i)},\mathbf{y}_k^{a,(i)}\}_{k=1}^{T}$ with NNs by minimizing the loss function:
		\[
		\mathcal{L}_F=\frac{1}{MT}\sum_{i=1}^{M}\sum_{k=1}^{T}\Vert x_k^{(i)} -F(\mathbf{y}_k^{a,(i)},\theta) \Vert _2 ^2.
		\]
		\label{algorithm1}
	\end{algorithm}

	\subsection{Online estimation}
	The surrogate model can be utilized for real-time estimation. Given the proxy transition model $K^{y,(t)}$, we compute  the posterior distribution $p(\hat{y}_k|Y_{k})$ with adaptive EnKF. Then the samples of $\mathbf{y}_k$ are collected with the trajectories $\{\mathbf{y}_k^{a,1},\cdots, \mathbf{y}_k^{a,N}\}$ defined in $(\ref{trajectory})$. We can get  the samples $\{x_k^{a,1},\cdots,x_k^{a,N}\}$ from the corresponding samples  $\{\mathbf{y}_k^{a,1},\cdots,\mathbf{y}_k^{a,N}\}$ through the reconstruction map $F$, and approximate the distribution $p(x_k|Y_k)$ by
	\[
	  p(x_k|Y_k)\approx  \frac{1}{N} \sum_{j=1}^{N} \delta_{x_k^{a,j}} (x_k).
	\]

	The prediction accuracy of $K^{y,(t)}$ may decrease for long-time prediction. To address this challenge, we update the DMD matrix with new observation in a new time window. Assuming that the DMD matrix varies with time, denoted by $K_k$, the initial $K_0$ is assigned as $K^{(t)}$. The filtering distribution of $p(K_k|Y_k)$ satisfies
	\[
	\begin{aligned}
		p(K_k|Y_k) &= \int p(K_k,\hat{y}_{1:k}|Y_k) d \hat{y}_{1:k}\\
		&= \int p(K_k|\hat{y}_{1:k}) p(\hat{y}_{1:k}|y_{1:k}) d \hat{y}_{1:k}\\
		&\approx p(K_k|\hat{y}_{1:k}^a),
	\end{aligned}
	\]
	where  the smooth distribution $p(\hat{y}_{1:k}|y_{1:k})$ is approximated by the posterior mean $\{\hat{y}_1^a,\cdots,\hat{y}_k^a\}$. We set the fixed-length time window $\{\hat{y}_{k-T+1},\cdots,\hat{y}_k\}$ to update $K_k$, i.e.,
	\[
	p(K_k|y_{k-T+1:k}) \approx p(K_k|\hat{y}_{k-T+1:k}^a),
	\]
	where $T$ represents the length of time window. 
	
	According to the definition of $\Psi(\hat{y}_k)$, we have
	$$
	\Psi(\hat{y}_{k+1})=K \Psi(\hat{y}_{k})+\hat{\Omega}_k,
	$$
	where $\hat{\Omega}_k \sim \mathcal{N}(0,\hat{M}_k)$ represents the model error between $\Psi(\hat{y}_{k+1})$ and $\Psi(\hat{y}_{k})$.
	We denote
	\[
	\begin{aligned}
		Y_{k-T+1:k-1}^{a}&:=[\hat{y}_{k-T+1}^a,\cdots,\hat{y}_{k-1}^a],\\
		Y_{k-T+2:k}^{a}&:=[\hat{y}_{k-T+2}^a,\cdots,\hat{y}_{k}^a],\\
		\hat{\Omega}_{k-T+1:k-1}&:=[\hat{\Omega}_{k-T+1},\cdots,\hat{\Omega}_{k-1}],\\
	\end{aligned}
	\]
	and then get
	\[
	\Psi(Y_{k-T+2:k}^{a})=K_k \Psi(Y_{k-T+1:k-1}^a) +\Omega_{k-T+1:k-1}.
	\]
	 We aim to obtain the least-squares estimate of $K_k$ by minimizing the following problem,
	\[
	K_k=\arg \min_{\tilde{K}} \Vert \Psi(Y_{k-T+2:k}^a)-\tilde{K} \Psi(Y_{k-T+1:k-1}^a) \Vert _{F}.
	\]
	Then the DMD matrix $K_k$ is updated using the formula,
	 \[
     K_k=\Psi(Y_{k-T+2:k}^{a}) \Psi(Y_{k-T+1:k-1}^a)^{\dagger}.
     \]
     Compared with the previous section, we update the DMD matrix in temporal space instead of the physical space to keep high accuracy for long-time prediction. The main steps of this algorithm are shown in Algorithm $\ref{algorithm2}$.
	
		\begin{algorithm}[H]
		\caption{Online estimation with the proxy model}
		\textbf{Input}: observation $y_k$ and $K_0$\\
		\textbf{Output}: posterior $p(x_k|Y_k)$\\
		~1:~ Estimate  $p(\hat{y}_k|Y_{k},K_{k-1})$ with adaptive EnKF and collect the samples of $\mathbf{y}_k$.\\
		~2:~ Estimate the posterior $p(x_k|Y_k)$ with the reconstruction map $F$:\\
		$$
		p(x_k|Y_k)=\frac{1}{N} \sum_{j=1}^{N} \delta_{F(\mathbf{y}_k^j,\theta)} (x_k).
		$$
		~3:~ Update $K_k$ with $K_k=\Psi(Y_{k-T+2:k}^{a}) \Psi(Y_{k-T+1:k-1}^a)^{\dagger}$ at predetermined time intervals. \\
		\label{algorithm2}
	\end{algorithm}
	
	\subsection{From EM perspective}\label{em_section}
	In this section, we develop  an  iterative process to build  the transition model $K^{y,(t)}$ from the perspective of EM algorithm. Denoting time-varying noise statistics as $M_k$ and $R_k$, our objective is to determine the maximum likelihood estimation for $K^{y}$. Since $K^y$ is determined by $K$, we  concentrate on deriving the solution for
	$K$. The cost function can be derived by
	\[
	\mathcal{J}(K)=-\log p(y_{1:T}|K)=-\log \int p(y_{1:T},z_{1:T}|K) \mathrm{d} z_{1:T},
	\]
	where $z_k=\Psi(\hat{y}_k)$ is the latent variable, and $ z_{1:T}=\{\Psi(\hat{y}_k)\}_{k=1}^{T}$. The optimization of this loss function can be achieved with EM algorithm.
	
	The initial DMD matrix $K^{(0)}$ can be calculated with the noisy observations. In the E step, it is essential to evaluate the distribution $p(z_{1:T}|y_{1:T},K^{y,(t)})$. To improve the  computational efficiency, the filtering distribution is utilized to approximate the smooth distribution $p(z_{1:T}|y_{1:T},K^{(t)})$. Given the transition model $K^{y,(t)}$, the posterior mean $\hat{y}_k^a$ can be obtained using adaptive EnKF.
	
	In the M step, our goal is to minimize
	\[
	\mathcal{J}(K)=-\int p(z_{1:T}|y_{1:T},K^{(t)})\log p(y_{1:T},z_{1:T}|K) \mathrm{d} z_{1:T}.
	\]
	We approximate $p(z_{1:T}|y_{1:T},K^{(t)})$  with the posterior mean trajectory $\{\Psi(\hat{y}_k^a)\}_{k=1}^{T}$ and get
	\[
	\begin{aligned}
		\mathcal{J}(K)&\approx-\log p(y_{1:T},z_{1:T}^{a}|K),
	\end{aligned}
	\]
	where $z_{1:T}^{a}=\{\Psi(\hat{y}_k^a)\}_{k=1}^{T}$. According to the form of $(\ref{dmd-ssm1})$, we can define the dynamical system of $\Psi(\hat{y}_{k})$ by
	\begin{equation}
		\begin{cases}
			\Psi_{k+1}=K \Psi_k+\hat{\Omega}_k,\\
			y_{k+1}=\tilde{h}
			(\Psi_k)+\eta_{k+1},
		\end{cases}
		\label{dmd-ssm2}
	\end{equation}
	where $\Psi_{k}:=\Psi(\hat{y}_{k})$. The observation operator $\tilde{h}$ contains the first $m$ components of $\Psi_k$ and
	\[
	\tilde{h}(\Psi_k)=
	\begin{bmatrix}
		\Psi_{k,1}\\
		\vdots\\
		\Psi_{k,m}\\
	\end{bmatrix}.
	\]
	Then according to $(\ref{dmd-ssm2})$, we can transfer $\mathcal{J}(K)$ into
	$$
	\begin{aligned}
		\mathcal{J}(K)
		&\approx-\log p(y_{1:T},\Psi_{1:T}^{a}|K)\\
		&\approx-\sum_{k=2}^{T}\log p(\Psi_k^a|\Psi_{k-1}^a)-\sum_{k=1}^{T}\log p(y_k|\Psi_k^a)-\log p(\Psi_1^a)\\
		&\approx \frac{1}{2} \sum_{k=2}^{T} \{ \Vert \Psi_k^{a}-K \Psi_{k-1}^a \Vert_{\hat{M}_k^{-1}}^2+\log |\hat{M}_k|+N_k \log(2\pi)\}-\log p(\Psi_1^a)\\
		&+ \frac{1}{2} \sum_{k=1}^{T} \{ \Vert y_k-\tilde{h}(\Psi_k^a) \Vert_{R_k^{-1}}^2+\log |R_k|+ m\log(2\pi)\}\\
		&\approx \frac{1}{2} \sum_{k=2}^{T}  \Vert \Psi_k^{a}-K \Psi_{k-1}^a \Vert_{\hat{M}_k^{-1}}^2+\cdots.
	\end{aligned}
	$$
	where $\Psi_{k}^a:=\Psi(\hat{y}_{k}^a)$. Except for the first term, the rest terms  are not related to  $K$. We denote the data matrices as
	$$
	\begin{aligned}
		&Y_{0}^{a}:=[\hat{y}_{1}^{a},\cdots,\hat{y}_{T-1}^{a}],\\
		&Y_{1}^{a}:=[\hat{y}_{2}^{a},\cdots,\hat{y}_{T}^{a}].\\
	\end{aligned}
	$$
	If we approximate the norm $\Vert \Psi_k^{a}-K \Psi_{k-1}^a \Vert_{\hat{M}_k^{-1}}$ with the Frobenius norm $\Vert \Psi_k^{a}-K \Psi_{k-1}^a \Vert_{F}$, the cost function can be approximated as follows:
	$$
	\mathcal{J}(K)\approx \min _{K} \Vert \Psi(Y_1^a)- K \Psi(Y_0^a)\Vert_{F}^2.
	$$
	The DMD matrix $K$ can be updated with $K^{(t+1)}=\Psi(Y_{1}^{a}) \Psi(Y_{0}^{a})^{\dagger}$.

	\section{KNN-Takens based EnKF}\label{section4}
	In this section, we present a nonparametric approach, namely KNN-Takens based EnKF, which employs analog forecasting strategies for online estimation. Compared with DMD-T EnKF, the transition model $f$ of $\hat{y}_k$ and the reconstruction map $F$ are constructed with nonparametric time series prediction
	methods for chaotic time series problems\cite{nonparametric,nonparametric-2,takens-filter}.
	
	In the offline stage, we consider constructing the proxy model for $f$ and denoise the observations. Based on the existing diffeomorphism from Takens' theorem, we embed the partial observation and utilize K-nearest neighbors (KNN) method for prediction. Given the posterior sample $\mathbf{y}_k^{a,j}=\{\hat{y}_{k}^{a,j},\cdots,\hat{y}_{k-d+1}^{a,j}\}$, we locate its $M$ nearest neighbors
	$$
	\{y_{k^{'}},\cdots,y_{k^{'}-d+1}\},\{y_{k^{''}},\cdots,y_{k^{''}-d+1}\},\cdots,\{y_{k^{M}},\cdots,y_{k^{M}-d+1}\}
	$$
	within the library of historical data with respect to Euclidean distance. The known successors $\{y_{k^{'}+1},y_{k^{''}+1},\cdots,y_{k^{M}+1}\}$ are used to predict the prior sample $\hat{y}_{k+1}^{f,j}$. Several analog forecasting operators are available, such as the locally constant analog operator (L-C operator) and the locally weighted linear regression operator (L-L operator)\cite{takens-filter,analog-EnKF}. These operators are defined as follows:
	
		\begin{itemize}
		\item L-C operator:
		$$
		f_{\text{LC}}(\mathbf{y}_k^{a,j})= w_1 y_{k^{'}+1} + w_2 y_{k^{''}+1}+\cdots+w_M y_{k^{M}+1},
		$$
		where $w_1,w_2,\cdots,w_M$ are weights that determine the contribution of each neighbor in constructing the prediction. The simplest form of weights is given by $w_1=w_2=\cdots=w_M=\frac{1}{M}$ and can be further defined as the normalized kernel function:
		$$
		w_i[\mathbf{y}_k^{a,j}]=\frac{g(\mathbf{y}_k^{a,j},\mathbf{y}^{(i)})}{\sum_{i=1}^{M}g(\mathbf{y}_k^{a,j},\mathbf{y}^{(i)})},
		$$
		where $\mathbf{y}^{(i)}$ represents the $i$-th nearest neighbor of $\mathbf{y}_k^{a,j}$ with $g(u,v)=\exp(-\lambda \Vert u-v \Vert^2)$.
		\item L-L operator: fit a weighted linear regression between  $\{\mathbf{y}_{k^{'}}^{a,j},\mathbf{y}_{k^{''}}^{a,j},\cdots,\mathbf{y}_{k^{M}}^{a,j}\}$ and the corresponding successors $\{y_{k^{'}+1},y_{k^{''}+1},\cdots,y_{k^{M}+1}\}$, where the weights are defined as the above content. We can calculate the regression slope $\theta(\mathbf{y}_{k}^{a,j})$ and intercept $\phi(\mathbf{y}_{k}^{a,j})$ using
		weighted least squares estimates. The L-L operator is written as
		$$
		f_{\text{LL}}(\mathbf{y}_k^{a,j})=\theta(\mathbf{y}_{k}^{a,j})\mathbf{y}_{k}^{a,j}+\phi(\mathbf{y}_{k}^{a,j}).
		$$
	\end{itemize}
	With the analog forecasting operators, the dynamical system can be approximated as follows:
	$$
	\begin{cases}
		\hat{y}_{k+1}=f_{\text{LC/LL}} (\mathbf{y}_{k})+\hat{\omega}_k,\\
		y_{k+1}=\hat{y}_{k+1}+\eta_{k+1}.
	\end{cases}
	$$
	Then the posterior distribution $p(\hat{y}_{k+1}|Y_{k+1})$ can be estimated with adaptive EnKF. We construct the set $\mathcal{Y}=\{\mathbf{y}_d^{a,(i)},\cdots,\mathbf{y}_T^{a,(i)}\}_{i=1}^M$, which comprises the posterior mean of $\mathbf{y}_k$ from the training labeled dataset. The ensemble of corresponding model state is assigned as $\mathcal{X}=\{x_d^{(i)},\cdots,x_T^{(i)}\}_{i=1}^M$. The two sets will be utilized to learn the reconstruction map online with analog forecasting methods.

	In the online stage, for a new observation $y_k$, we firstly estimate $p(y_{k}|Y_{k})$ with adaptive EnKF. Then the samples of $\mathbf{y}_k$ are collected with $\{\mathbf{y}_k^{a,1},\cdots, \mathbf{y}_k^{a,N}\}$ defined in $(\ref{trajectory})$. Given $\mathbf{y}_{k}^{a,j}$, we locate its $M$ nearest neighbors within the set $\mathcal{Y}$:
	$$
	\{\mathbf{y}_{k^{'}}^{a,j}, \mathbf{y}_{k^{''}}^{a,j}, \cdots,\mathbf{y}_{k^{M}}^{a,j}\}.
	$$
	The corresponding labels of these $M$ neighbors in the set $\mathcal{X}$ are  $\{x_{k^{\prime}}, x_{k^{\prime\prime}},\cdots, x_{k^M}\}$. We can utilize L-C operator or L-L operator to predict $x_k^{a,j}$, respectively. The selection of these operators mainly depends on the available computational resources and the complexity of the application\cite{analog-EnKF}. The framework of KNN-Takens based EnKF is shown in Algorithm $\ref{algorithm3}$.
	
	\begin{algorithm}[H]
		\caption{KNN-Takens based EnKF}
		\textbf{Input}: historical trajectory $\{x_k^{(i)},y_k^{(i)}\}_{k=1}^T, i=1,2,\cdots,M$, new observation $y_n$\\
		\textbf{Output}: posterior $p(x_{n}|Y_{n})$\\
		~1:~ In the offline stage, denoise the observations with analog forecasting operators and construct the sets
		$$
		\begin{aligned}
			\mathcal{Y}&=\{\mathbf{y}_d^{a,(1)},\cdots,\mathbf{y}_T^{a,(1)},\cdots,\mathbf{y}_d^{a,(M)},\cdots,\mathbf{y}_T^{a,(M)}\},\\
			\mathcal{X}&=\{x_d^{(1)},\cdots,x_T^{(1)},\cdots,x_d^{(M)},\cdots,x_T^{(M)}\}.
		\end{aligned}
		$$  \\
		~2:~ In the online stage, estimate $p(y_{n}|Y_{n})$ with adaptive EnKF and estimate $p(\mathbf{y}_{n}|Y_{n})$.\\
		~3:~ For each sample $\mathbf{y}_{n}^{a,j}$, locate its $M$ nearest neighbors and utilize the analog forecasting operator to estimate $x_{n}^{a,j}$. The posterior distribution of $x_{n}$ is approximated by
		$$
		p(x_{n}|Y_{n})\approx \frac{1}{N} \sum_{j=1}^{N} \delta_{x_{n}^{a,j}} (x_{n}).
		$$
		\\
		\label{algorithm3}
	\end{algorithm}

\section{Numerical results}\label{section5}
In this section, we implement some numerical examples to illustrate the performance  of DMD-T EnKF and KNN-T EnKF. In Section \ref{sec5.1}, we consider a pendulum model and estimate the overall model state from the partial and noisy observation by DMD-T EnKF. In Section \ref{sec5.2}, we consider a stochastic model for topographic mean flow. In Section \ref{sec5.3}, KNN-T EnKF is applied to a Lorenz 63 model. In Section \ref{sec5.4}, we consider an example for Allen-Cahn equation with partial observations using DMD-T EnKF. To evaluate the performance of online estimation, we utilize the following metrics:
$$
\begin{aligned}
	\text{RMSE}&=\sqrt{\frac{1}{nT}\sum_{k=1}^{T}\Vert x_k-x_k^a\Vert_2^2},\\
	\text{Spread}&=\sqrt{\frac{1}{nT}\sum_{k=1}^{T}\text{Tr}(P_k^a)},
\end{aligned}
$$
where $k$ denotes the time index and $n$ is the dimensionality of $x_k$. $x_k^a$ and $P_k^a$ represents the posterior mean and covariance of the model state, respectively.

\subsection{Pendulum Model}\label{sec5.1}
	\begin{figure}[h]
		\centering
		\includegraphics[scale=0.55]{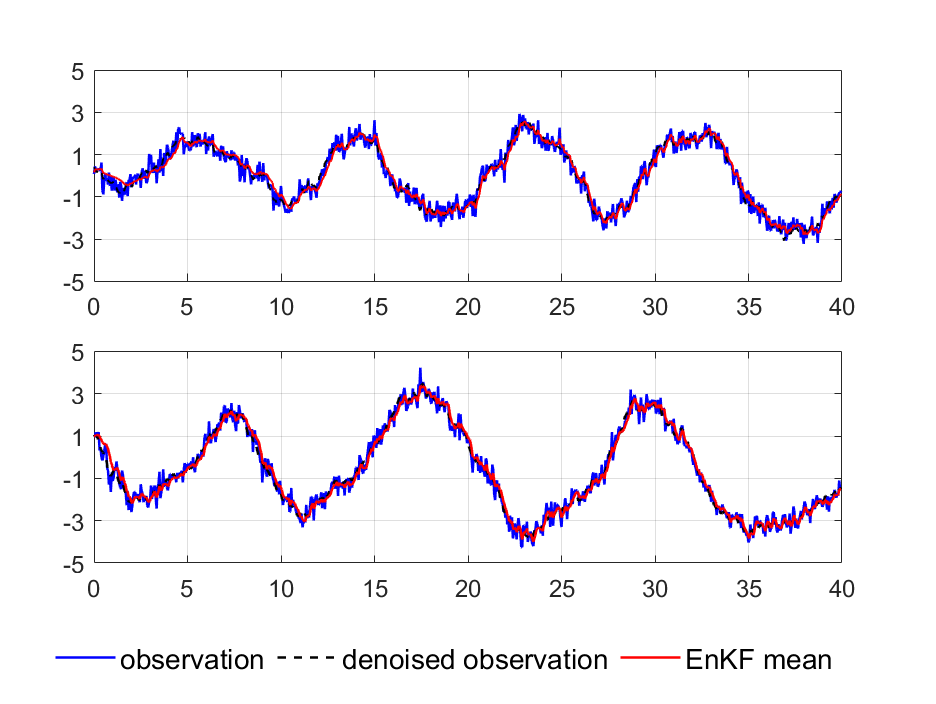}
		\caption{The posterior mean $\hat{y}_k^a$ (red), the denoised observation (black dashed), the observation (blue) for the historical trajectories.  We calculate $K^{y,(0)}$ with noisy snapshots and update the matrix $K^{y,(1)}$ with the posterior mean $\hat{y}_k^a$. The $l^2$ norm of $K^{y,(0)}-K^{y}$ and $K^{y,(1)}-K^{y}$ are 0.1058 and 0.0126, respectively. }
		\label{danbai_1}
	\end{figure}
	In this subsection, we consider the pendulum model with  the state vector $x_k=[\omega_k,\phi_k]^{T}$, which incorporates the angular velocity $\omega_k$ and the angle $\phi_k$.  Suppose  that the evolution process is influenced by a random  noise, then the  transition model are governed by the following equation,
	\[
	\left\{
	\begin{aligned}
		\frac{\partial \omega}{\partial t} &= -\frac{g}{l} \sin (\phi)+\sigma_1 \frac{dW_1}{dt},\\
		\frac{\partial \phi}{\partial t} &= \omega+\sigma_2 \frac{dW_2}{dt},
	\end{aligned}
	\right.
	\]
	where $g$ denotes the gravitational acceleration with a constant value of $9.8$  $[m/sec^{2}]$, and $l$ represents the length of the pendulum's string and is set as $l=20$. The continuous equation is discretized into a discrete-time system using Euler-Maruyama's method. The model can be expressed as
	\[
	\left\{
	\begin{aligned}
		\omega_{k+1}&=\omega_k -\Delta_t \frac{g}{l} \sin (\phi_k)+\sigma_1 \xi_{1,k}\\
		\phi_{k+1}&=\phi_k+\Delta_t \omega_k+\sigma_2 \xi_{2,k},
	\end{aligned}
	\right.
	\]
	where $\Delta t=0.05$ and $\xi_{1,k}$,  $\xi_{2,k}\sim N(0,1)$, with noise variances $\sigma_1^2=\sigma_2^2=0.002$. Additionally, the observation function is defined as
	$$
	y_k=\phi_k + \eta_k,
	$$
	where $\eta_k \sim N(0,R)$ and $R=0.1$. The historical trajectories consist of $\{x_k^{(i)},y_k^{(i)}\}_{k=1}^{800},\ \ i=1,2,\cdots,40$. In the offline stage, the surrogate model $K^{y}$ and the reconstruction map $F$ is constructed by Algorithm $\ref{algorithm1}$. We calculate $K^{y,(0)}$ with noisy snapshots and update the matrix $K^{y,(1)}$ with the posterior mean $\hat{y}_k^a$. The $l^2$ norm of $K^{y,(0)}-K^{y}$ and $K^{y,(1)}-K^{y}$ are 0.1058 and 0.0126, respectively. We select two historical trajectories.  The posterior mean $\hat{y}_k^a$ is showed  in Fig.\;$\ref{danbai_1}$. It illustrates that the posterior mean $\hat{y}_k^a$ can track the truth $\hat{y}_k$, achieving noise reduction. Additionally, the posterior mean $\hat{y}_k^a$ can provide more accurate estimates of the transition model. After that, we set the delay coordinate vector $\mathbf{y}_k$ with different lengths $d=2,5,10,15$. The reconstruction map is learned using the dataset  $\{x_k^{(i)},\mathbf{y}_k^{(i),a}\}_{k=1}^{800}$.
	
\begin{figure}[htp]
		\centering
		\hspace{-10pt}
		\subfigure[]{
			\includegraphics[height=8.5cm, width=9.0cm]{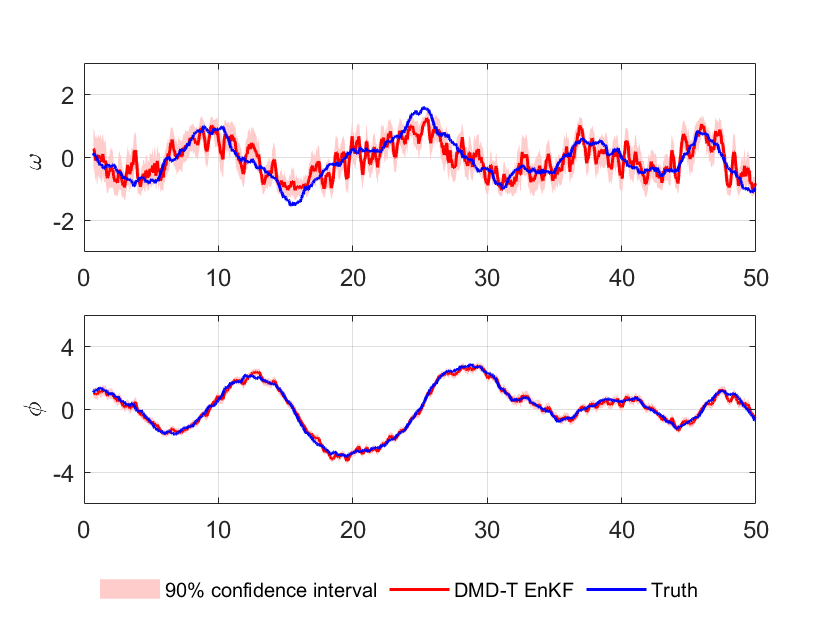}
		}
		\subfigure[]{
			\includegraphics[height=8.5cm, width=7cm]{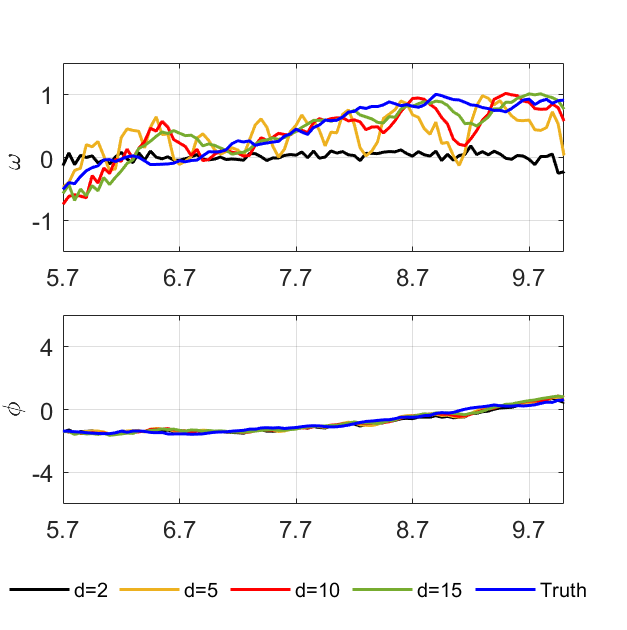}
		}
		
		\vspace{-5pt}
		\caption{(a) The posterior mean $x_k^a$ of test data by DMD-T EnKF when $d=10$. (b) The posterior mean  $x_k^a$ of test data for different lengths $d=2,5,10,15$ in the time window $[5.7,9.7]$.}
		\label{danbai_2}
	\end{figure}

	In the online stage, we utilize test data comprised of $\{y_k\}_{k=1}^{1000}$ with different initial value to test the proxy model. The posterior mean $x_k^a$ of test data is illustrated  in Fig.\;$\ref{danbai_2}$. Setting $d=10$, the posterior mean $x_k^a$ is able to track the true trajectory and the truth falls within the confidence interval. The estimation of $\phi_k$ is more accurate and exhibits less uncertainty than that of $\omega_k$. This discrepancy is attributed to the observation operator only containing information related to $\phi_k$. From Fig.\;$\ref{danbai_2}$(b), as $d$ increases, the posterior mean $\omega_k^a$ gradually approaches the truth value.
	
\begin{figure}[htp]
		\hspace{-10pt}
		\includegraphics[scale=0.5]{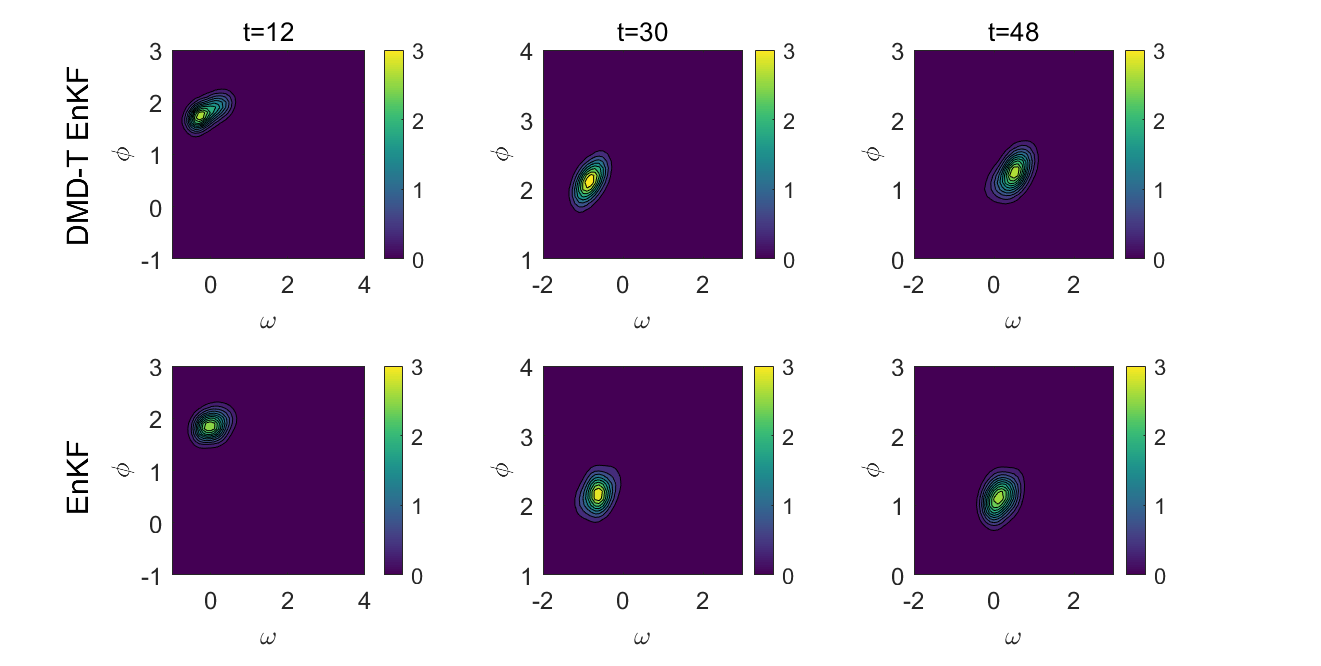}
		\caption{The filtering distribution $p(x_k|Y_k)$ of test data at time points $t=12,30,48$ compared to EnKF ($d=10$). The horizontal axis represents $\omega$, and the vertical axis represents $\phi$ estimated through KDE with the same bandwidth in two methods.}
		\label{danbai_3}
	\end{figure}

     Fig.\;$\ref{danbai_3}$  shows  the filtering distribution $p(x_k|Y_k)$ of test data through kernel density estimation (KDE), where the horizontal axis represents $\omega_k$ and the vertical axis represents $\phi_k$. It can be observed that the MAP estimate by DMD-T based EnKF is consistent with that of EnKF. The posterior distributions by DMD-T EnKF also exhibit a similar pattern  to those of EnKF.
	
	\begin{figure}[htp]
		\centering
		\includegraphics[scale=0.45]{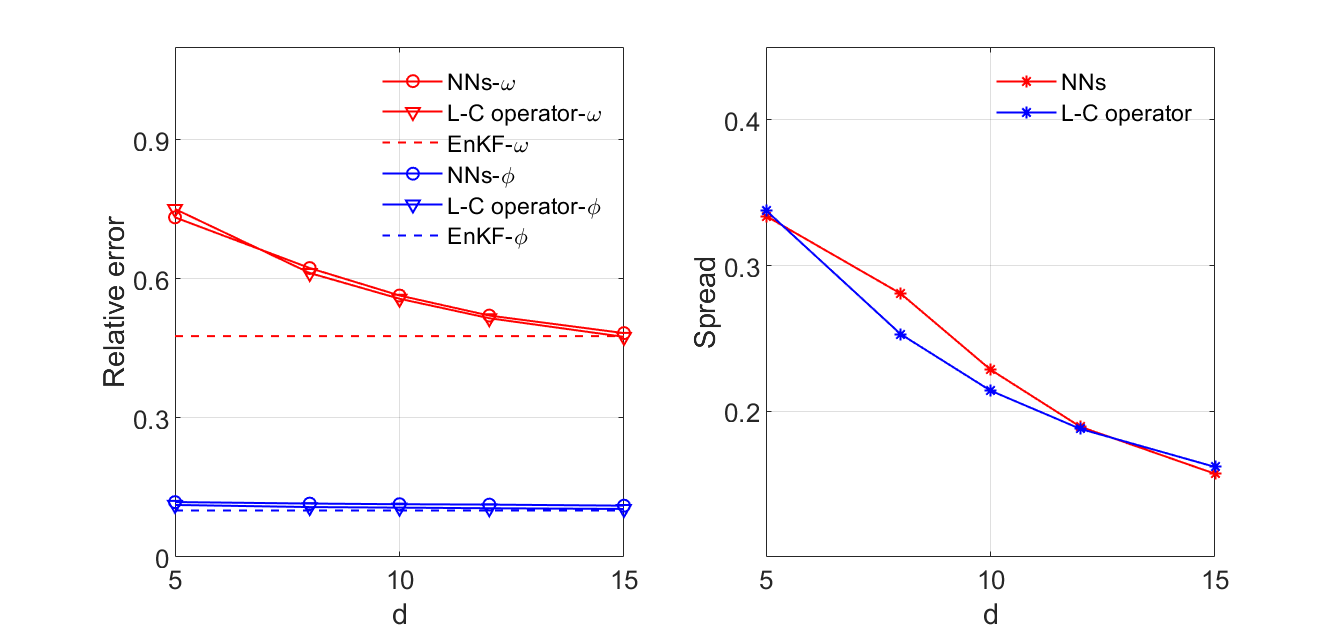}
		\caption{The relative errors and spreads with different modeling methods for the reconstruction map: NNs, L-C operator.}
		\label{danbai_4}
	\end{figure}
	To explore the impact of the time-delay  length $d$, we compute  the relative errors and spreads of NNs and L-C operator.     We plot the errors in Fig.\;$\ref{danbai_4}$ and see that the relative errors decrease  as $d$ increases. The relative error of $\phi_k$ is not strongly correlated with the length $d$  because the current observation has contained the complete information of $\phi_k$. While the relative error of $\omega_k$ exhibits highly correlation with $d$ and an increase of $d$ leads to better  estimation accuracy. Furthermore, the spreads decrease with increase of $d$. This may be due  to the increased input dimension, which makes the reconstruction map more robust. With different modeling methods for the reconstruction map, the difference of NNs and L-C operator is not  significant. When the model state dimension is low, nonparametric methods can be utilized to estimate $x_k$ without the explicit learning of the map $F$.

	\subsection{ A stochastic model}\label{sec5.2}
			In this subsection, we consider a stringent paradigm model for topographic mean flow interaction, which is a kind of physics constrained multi-level regression models\cite{sde-1,sde-2}. With a specific choice of coefficients, the system can be expressed as follows:
	\begin{equation}
	dx=(Lx-Dx+B(x,x))dt+\sigma \Gamma dW,
	\label{sde}
	\end{equation}
	where $x=(u,v,w)^{T},B(x,x)=(0,auw,-auv)^{T},dW=(dW_1,dW_2)^{T}$,
	$$
	L=\left(
	\begin{array}{lll}
		0 & \omega & 0\\
		-2\omega &0 &-\beta\\
		0 & \beta & 0
	\end{array}
	\right), \quad
	D=\left(
	\begin{array}{lll}
		0 & 0 & 0\\
		0 &\gamma &0\\
		0 & 0 & \gamma
	\end{array}
	\right), \quad
	\Gamma=\left(
	\begin{array}{ll}
		0 & 0\\
		1&0 \\
		0&1
	\end{array}
	\right).
	$$
	This form is similar to the Charney and DeVore model for nonlinear regime behavior without dissipation and forcing\cite{sde-1,sde-2}. We numerically discretize the continuous model with the Euler-Maruyama scheme as follows:
	$$
	x_{i+1}=x_i+\delta t (Lx_i-Dx_i+B(x_i,x_i))+\sigma \sqrt{\delta t} \Gamma \Delta W_{i+1}.
	$$
	where $\Delta W_{i+1}\sim N(0,I_2)$. The parameters are assigned as $\omega=3/4,\gamma=1/2,\beta=1,a=1,\sigma=1/\sqrt{2}$ and $\delta t=0.1$. Assume the noisy observation are available from the observation function
	$$
	y_i=h(x_i)+\eta_i,
	$$
	where $\eta_i\sim N(0,R)$. The functional form of $h$ is specified as:
	$$
	h=\left(
	\begin{array}{l}
		u+v\\
		2w\\
	\end{array}\right),\quad
	R=\left(
	\begin{array}{ll}
		r  & 0\\
		0 &r
	\end{array}\right),
	$$
	where $r=0.1$.
	
	\begin{figure}[htp]
		\centering
		\hspace{-50pt}
		\subfigure[]{
			\includegraphics[height=7.5cm, width=10cm]{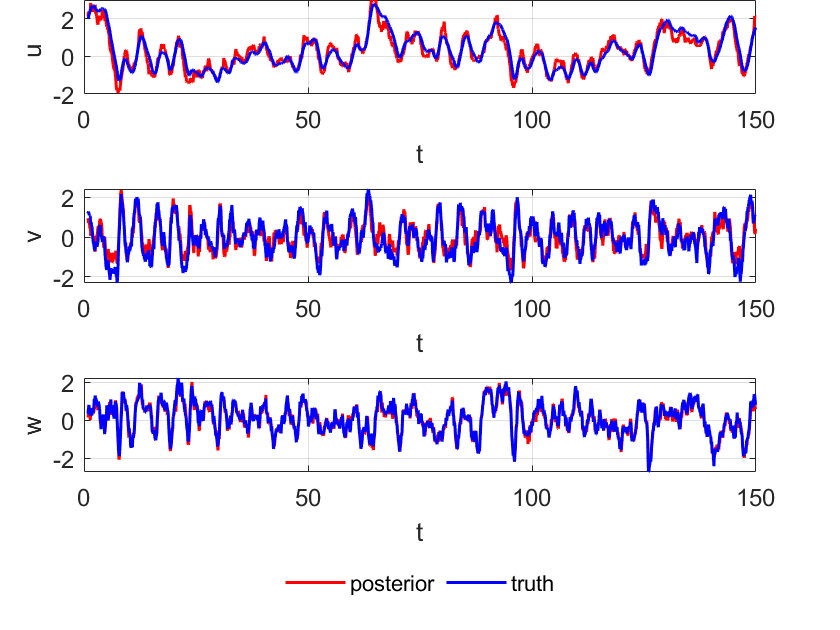}
		}
		\subfigure[]{
			\includegraphics[height=8.5cm, width=6cm]{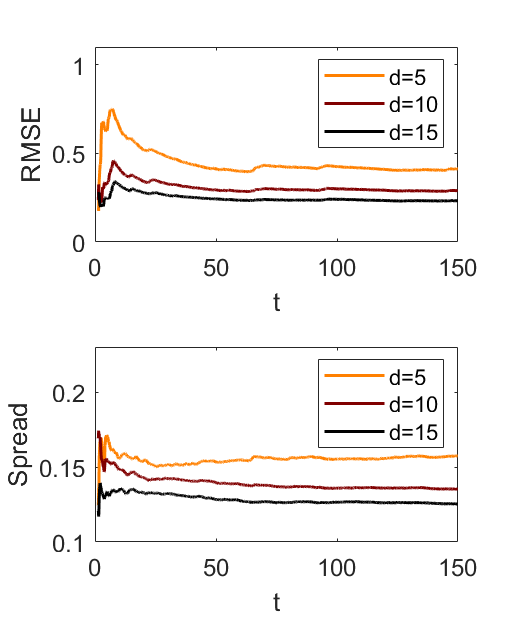}
		}
		
		\vspace{-5pt}
		\caption{(a) The posterior mean $x_k^a$ of test data with DMD-T EnKF when $d=10$. (b) RMSEs and spreads of test data with different lengths $d=5,10,15$.}
		\label{sde_2}
	\end{figure}
	
	\begin{figure}[h]
		\centering
		\includegraphics[scale=0.5]{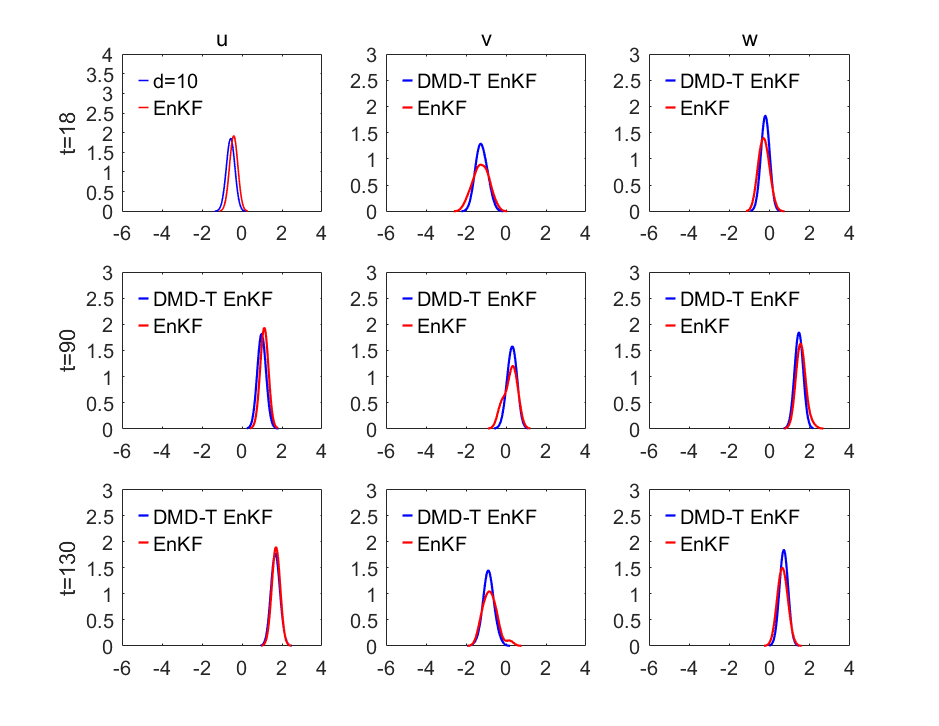}
		\caption{The marginal distribution of $p(x_k|Y_k)$ estimated at time points $t=18,90,130$ compared to EnKF ($d=10$). The horizontal axis represents the value of variables, and the vertical axis represents the probability density estimated through KDE with the same bandwidth in two methods.}
		\label{sde_3}
	\end{figure}

	\begin{table}[h]
		\centering
		\begin{tabular}{c|c|c|c|c}
			\hline
			& $K^{y,(0)}$ & $K^{y,(1)}$ & \textit{$K^{y,(2)}$} & $K^{y,(3)}$  \\
			\hline
			error($l^2$) & 0.059999  & 0.018666 & 0.015302 & 0.015457 \\
			\hline
		\end{tabular}
		\caption{the $l^2$ error between $K^{y,(t)}$ and $K^{y}$, $t=0,1,2,3$. }
		\label{sde_1}
	\end{table}

	Without the underlying dynamics, our aim is to construct the surrogate model for online estimation. The historical trajectories is composed of $\{x_k^{(i)},y_k^{(i)}\}_{k=1}^{500},\ \ i=1,\cdots,20$. Firstly, we construct the surrogate model with Algorithm $\ref{algorithm1}$. We repeat EnKF-DMD three times to estimate the DMD matrix $K^{y}$ and the denoised observation $\hat{y}$. The $l^2$ error of $K^{y,(t)}-K^{y}$ is in Table $\ref{sde_1}$. We observe that the iteration process can effectively reduce the $l^{2}$ error between the estimated value and the truth value. This improvement can be explained  from the perspective of the EM framework.
	
	Then we set the delay coordinate vector $\mathbf{y}_k=\{\hat{y}_k,\cdots,\hat{y}_{k-d+1}\}$ with different time-delay  lengths $d=5,10,15$.  We test the proxy model with new observation $\{y_1,y_2,\cdots,y_{1500}\}$ generated from different initial value. The posterior mean $x_k^a$ of the test data is displayed in Fig.\;$\ref{sde_2}$. The posterior mean can accurately track the truth. We observe that the estimation of $w$ achieves better accuracy than  other variables, due to the observation function containing complete information of $w$. Additionally, RMSEs and spreads decrease with increase of $d$.   We plot  marginal distribution of the three variables at $t=15,90,130$ in Fig.\;$\ref{sde_3}$. The horizontal axis represents the value of variables, and the vertical axis represents the probability density estimated through KDE. It can be observed that the distribution estimated by DMD-T EnKF is very close to the distribution by  EnKF.
	
	 To explore the sensitivity with different noise levels, we compute RMSEs and spreads compared to EnKF for $\gamma=0.1, 0.15, 0.2, 0.25, 0.3$. The result  is listed  in Table $\ref{sde_4}$. With the increase of noise magnitude, the RMSE increase rate of DMD-T EnKF remains approximately  $4\%$, but  the standard EnKF  has  $10\%$. This indicates that our proposed method remains effective under high noise magnitude conditions. Additionally, the spreads by DMD-T EnKF are smaller than those by EnKF. Unlike the standard EnKF, which provides prior information in the state space, DMD-T EnKF provides the prior information in the denoised observation space. The uncertainty of $x_k$ depends on the uncertainty of $\mathbf{y}_k$. The decrease in uncertainty may be attributed to the robustness of the reconstruction map, which is trained with NNs.
	\begin{table}[h]
		\centering
		\begin{tabular}{c|cc|cc}
			\hline
			& \multicolumn{2}{c|}{DMD-T EnKF}          & \multicolumn{2}{c}{EnKF}                \\ \hline
			noise magnitude & \multicolumn{1}{c|}{RMSE}     & Spread   & \multicolumn{1}{c|}{RMSE}     & Spread   \\ \hline
			0.1             & \multicolumn{1}{c|}{0.291055} & 0.135308 & \multicolumn{1}{c|}{0.184065} & 0.17918  \\
			0.15            & \multicolumn{1}{c|}{0.303017} & 0.156167 & \multicolumn{1}{c|}{0.214888} & 0.209518 \\
			0.2             & \multicolumn{1}{c|}{0.319642} & 0.176229 & \multicolumn{1}{c|}{0.238832} & 0.233021 \\
			0.25            & \multicolumn{1}{c|}{0.330295} & 0.188462 & \multicolumn{1}{c|}{0.258612} & 0.252394 \\
			0.3             & \multicolumn{1}{c|}{0.335907} & 0.201737 & \multicolumn{1}{c|}{0.275561} & 0.268970 \\ \hline
		\end{tabular}
		\caption{RMSEs and spreads of test data with different noise magnitude $\gamma$ (d=10).}
		\label{sde_4}
	\end{table}

	\subsection{Lorenz 63}\label{sec5.3}
	In this section, we consider a chaotic system, which is explored by Lorenz for atmosphere convection in 1963. This  Lorenz model exhibits high sensitivity to initial conditions and has chaotic characteristics. The model is described by a coupled nonlinear ordinary differential equations with an additive noise. The solution, denoted as $v=(v_1,v_2,v_3)\in \mathbb{R}^3$, satisfies
	\[
	\left\{\begin{array}{l}
		\frac{d v_1}{d t}=a\left(v_2-v_1\right)+\gamma_1\frac{dW_1}{dt}, \\
		\frac{d v_2}{d t}=-a v_1-v_2-v_1 v_3+\gamma_2\frac{dW_2}{dt}, \\
		\frac{d v_3}{d t}=v_1 v_2-b v_3-b(r+a)+\gamma_3\frac{dW_3}{dt}. \\
	\end{array}\right.
	\]
	The parameter values are assigned as $(a,b,r)=(10,\frac{8}{3},28)$ and $\gamma_1=\gamma_2=\gamma_3=0.05$.
	
	We can solve  this continuous system with Euler-Maruyama scheme, where $\delta_t=0.01$. It is assumed that noisy observation can be obtained through an observation function, which satisfies
	$$
	y_k=h(v_k)+\eta_k,
	$$
	where $\eta_k \sim \mathcal{N}(0,\sigma^2), \sigma^2=2$. The historical trajectories is comprised of  $\{x_k^{(i)},y_k^{(i)}\}_{k=1}^{200},\ \ i=1,\cdots,30$. Algorithm $\ref{algorithm3}$ is employed for online estimation with nonparametric time series prediction methods. We locate $N=100$ nearest neighbors for analog forecasting. The delay coordinate vectors $\mathbf{y}_k$ is formed with the length $d=20$. The set of delay coordinate vectors $\mathcal{Y}$ with the ensemble of successors $\mathcal{X}$ is assigned as training data for the  nonparametric online estimation.
	\begin{figure}[!h]
		\centering
		\includegraphics[scale=0.5]{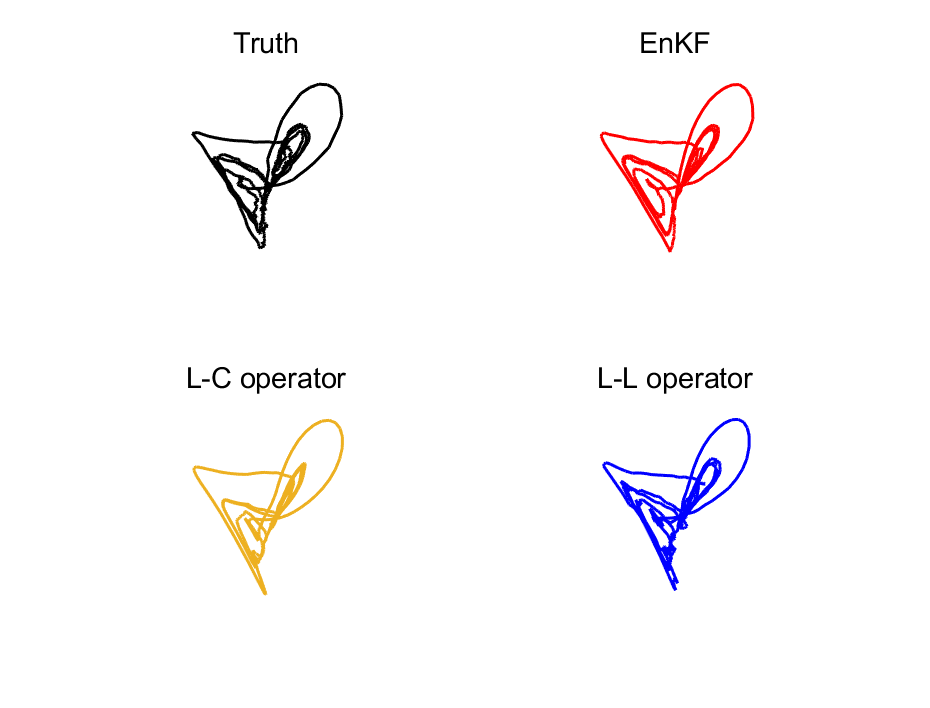}
		\caption{The posterior mean $v_k^a$ of test data with KNN-T EnKF. Two different operators are utilized for the transition model and the reconstruction map: L-C operator and L-L operator, respectively. We compute the two metrics with different operators: $\text{RMSE}_{\text{L-C}}=2.9591, \text{RMSE}_{\text{L-L}}=2.2205$ and $\text{spread}_{\text{L-C}}=0.3719, \text{spread}_{\text{L-L}}=1.0703$.}
		\label{lorenz63_1}
	\end{figure}
	
	\begin{figure}[!h]
		\centering
		\includegraphics[scale=0.4]{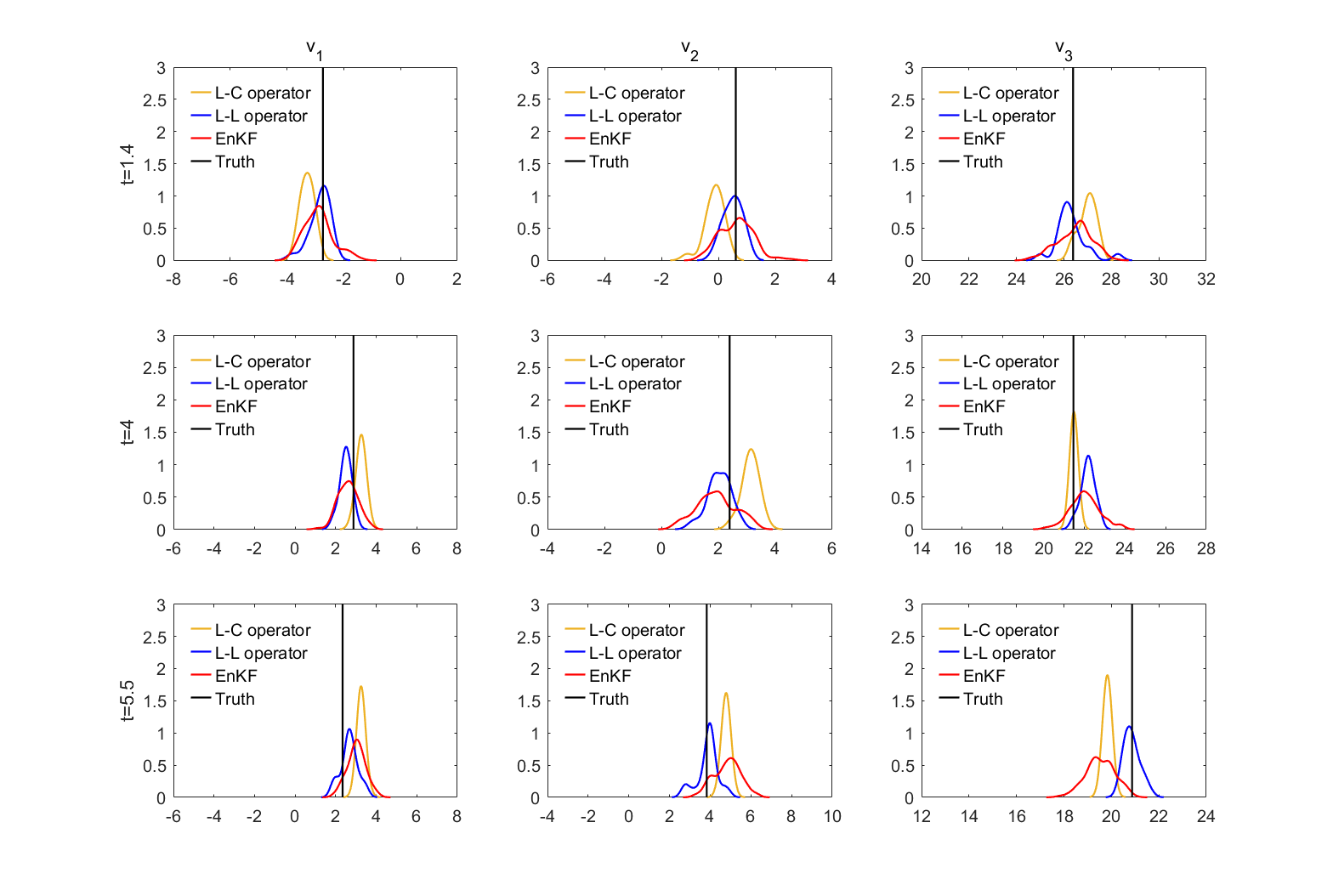}
		\caption{The marginal distribution of $p(v_k|Y_k)$ estimated at time points $t=1.4, 4, 5.5$ compared to EnKF ($d=10$). The horizontal axis represents the value of variables, and the vertical axis represents the probability density estimated through KDE with the same bandwidth among three methods.}
		\label{lorenz63_2}
	\end{figure}
	
	Suppose an observation operator is given by $h_1(v)=v_1+v_2+v_3$, which includes only partial information about $v$. The surrogate model is constructed with Algorithm $\ref{algorithm3}$ and we utilize new observation $\{y_k\}_{k=1}^{700}$ to test the proxy model. We estimate the filtering distribution $p(v_k|Y_k)$ using KNN-T EnKF based on  L-C operator and L-L operator, respectively. The result is shown in Fig.\;$\ref{lorenz63_1}$. It illustrates that our method can track the truth and effectively estimate the chaotic attractors. We can find that the estimation accuracy of $\text{L-L operator}$ is higher than L-C operator. The $\text{L-L operator}$ is semi-parametric with a locally weighed linear regression, while L-C operator is purely parametric relying on the training data. The metric spread measures the uncertainty inherent in samples $x_k^{a,i}$. The spread of L-L operator is larger than that of L-C operator. It leads to a big dispersion of prior ensembles by the posterior ensembles, with the transition model learned by L-L operator in the observation space.
	\begin{figure}[h]
		\centering
		\includegraphics[scale=0.5]{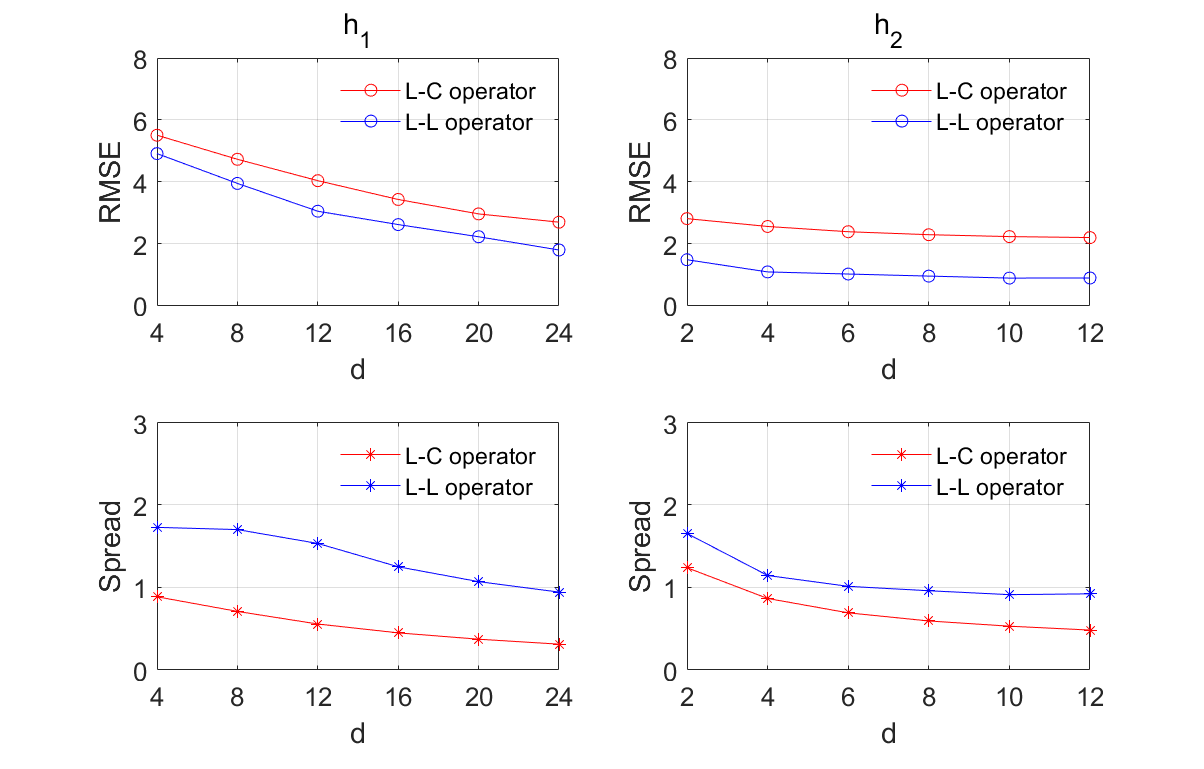}
		\caption{RMSEs and spreads with two analog forecasting operators for different observation functions $h_1$ and $h_2$.}
		\label{lorenz63_3}
	\end{figure}
	
	Fig.\;$\ref{lorenz63_2}$ shows a comparison of filtering distribution with KDE by L-C operator, L-L operator, EnKF and the ground truth when $t=1.4, 4, 5.5$. Compared to standard EnKF, the two operators can both provide relatively accurate MAP estimations that are close to the ground truth. Our method estimates the uncertainty of distribution in the delay observation space instead of the original space. The distribution estimated by L-L operator will be more dispersed than L-C operator with respect to the lager spread and sensitivity of linear regression. Due to the nonlinearity of map $F$, the distribution exhibits stronger non-Gaussian characteristics.
	
	To explore the impact of different observation functions, we suppose that new observation operator $h_2$ satisfies
	$$
	h_2(v)=\left(
	\begin{aligned}
		v_1+&v_2+v_3\\
		&v_3
	\end{aligned}\right).
	$$
	We calculate RMSEs and spreads for two analog forecasting operators with different observation function $h_1$ and $h_2$. The trend of RMSE and spread with regard to  $d$ is described in Fig.\;$\ref{lorenz63_3}$. The observation's dimension of $h_2$ is twice that of $h_1$. So we establish the time-delay length of $h_2$ to be half that of $h_1$, which can ensure the same dimension of $\mathbf{y}_k$ for two observation functions. From Fig.\;$\ref{lorenz63_3}$, RMSEs and spreads are expected to decrease as $d$ increases, attributed to the augmented information encapsulated within the delay coordinate vectors. If $d$ is too small, the information derived from $\mathbf{y}_k$ may be inadequate for estimating the chaotic attractors. Under a fixed dimension of $\mathbf{y}_k$, it becomes evident that the estimation of $h_2$ yields better accuracy  compared to $h_1$. It is rational to assert that $h_2$ contains information pertaining to
	$v_3$, whereas $h_1$ encompasses only partial information regarding the sum of three components in $v$. The estimation accuracy of L-L operator is larger than that of L-C operator and provides larger estimation uncertainty.
 	\begin{figure}[h]
 		\centering
 		\includegraphics[scale=0.45]{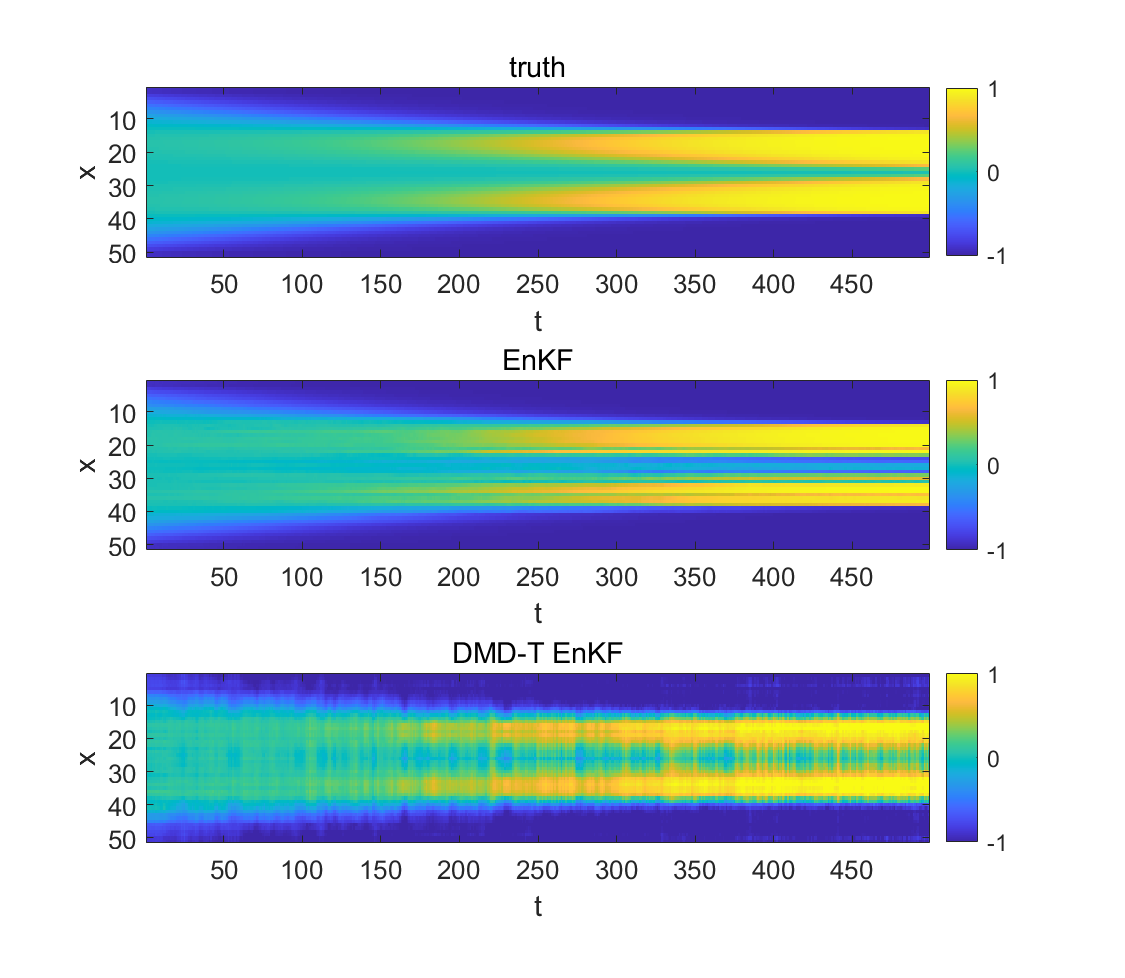}
 		\caption{The posterior mean $u_k^a$ of test data using DMD-T EnKF and EnKF.}
 		\label{ac_1}
 	\end{figure}

	\subsection{Allen-Cahn equation}\label{sec5.4}
	Allen-Cahn equation has been widely used for  phase transitions and interfacial dynamics in materials science. Initially proposed by Allen and Cahn, this equation was introduced to describe the motion of anti-phase boundaries within crystalline solids. It has found extensive application in various domains including crystal growth, phase transitions, image analysis, grain growth, and numerous other areas\cite{a-c,a-c-1}.
	
	In this section, we consider the one-dimensional Allen-Cahn equation
	\[
	\left\{\begin{array}{l}
		\frac{\partial u}{\partial t}=\epsilon \frac{\partial^2 u}{\partial x^2}-\theta (u^3-u);, \\
		u(-1,t)=u(1,t), \quad u_x(-1,t)=u_x(1,t)\\
		u(x,0)=g(x), \\
	\end{array}\right.
	\]
	where $x\in[-1,1],t\in[0,1]$. The parameters are set as $\epsilon=0.0001,\theta=5$.   We employ numerical discretization techniques, utilizing the temporal Crank-Nicolson method for time discretization and a second-order central difference scheme for spatial discretization. The spatial domain $[-1,1]$ is discretized using a finite grid with a spatial increment of $\Delta x=1/50$, while the temporal domain $[0,1]$ is discretized with a time step of $\Delta t=1/500$.  We get  a reference solution by employing the finite element method across all grid points.
	
	\begin{figure}[h]
		\centering
		\includegraphics[scale=0.4]{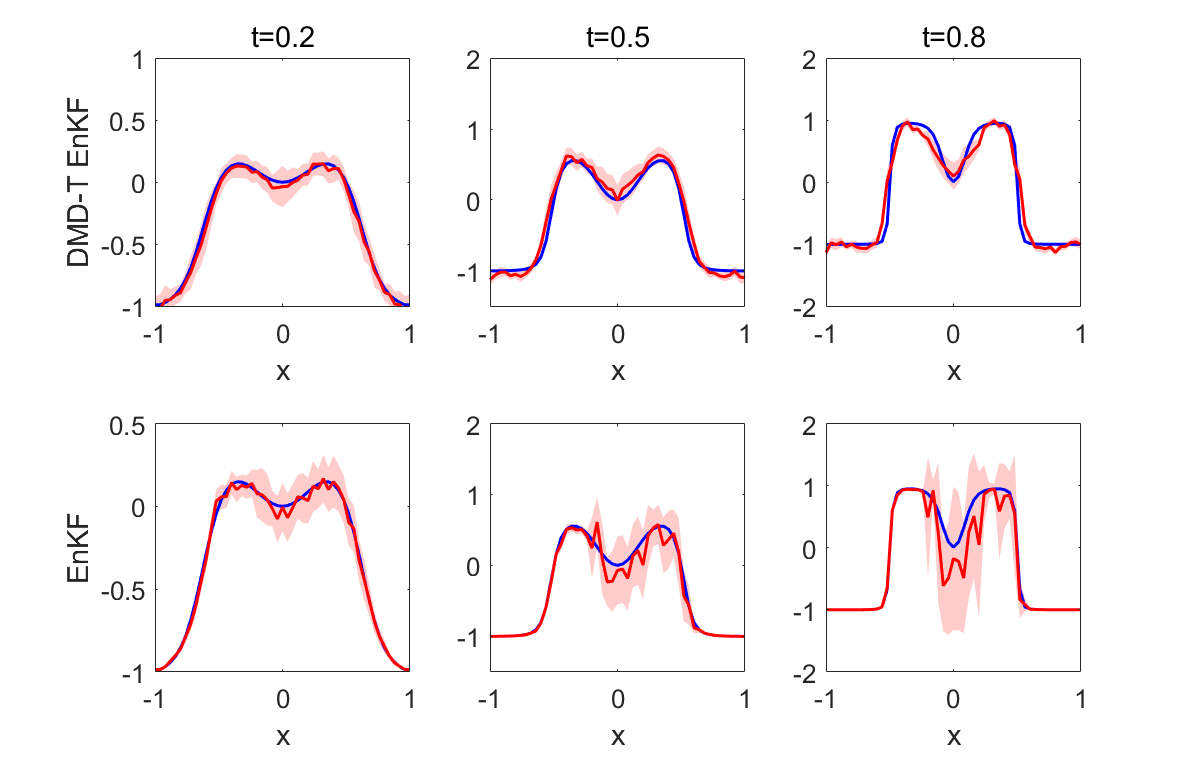}
		\caption{The filtering distribution of test data estimated through KDE at time points $t=0.2,0.5,0.8$.}
		\label{ac_2}
	\end{figure}
	\begin{figure}[h]
		\centering
		\includegraphics[scale=0.45]{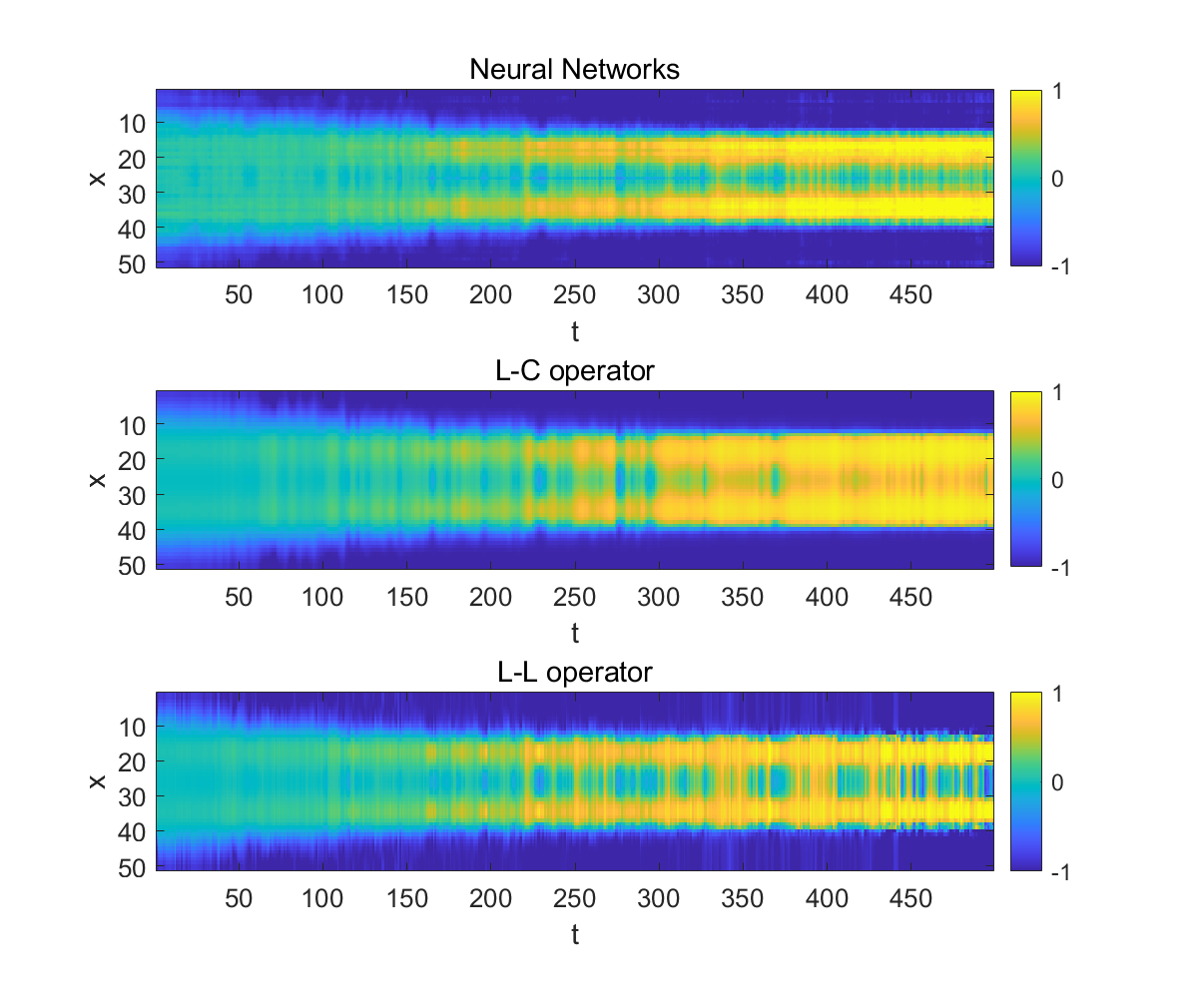}
		\caption{The posterior mean $u_k^a$ of test data estimated with different operators for the reconstruction map: NNs, L-C operator and L-L operator.}
		\label{ac_3}
	\end{figure}
	
	Consider that the observation of this equation is accessible through an observation function
	$$
	y=h(u)+\eta,
	$$
	where $\eta \sim \mathcal{N}(0,\sigma I), \sigma=0.1$. We assume $h$ contains only the first 20 components of $u$, which satisfies
	$$
	h(u)=u_{1:20}.
	$$
	The historical data consists of $\{u_k^{(i)}, y_k^{(i)}\}_{k=1}^{501},\ \ i=1,\cdots,10$ initialized with different initial conditions $u(x,0)$. Algorithm $\ref{algorithm1}$ is utilized to learn the proxy model with historical data. The time-delay length is $d=3$. We employ the test data $\{y_k\}_{k=1}^{501}$, derived from the solution with initial conditions $u(x,0)=x^2\cos(\pi x)$, to examine the model's performance. Fig.\;$\ref{ac_1}$ depicts the posterior mean $u_k^a$ for the test data using DMD-T based EnKF and EnKF. We observe that the model state $u_k$ estimated by the proposed method demonstrates symmetry about the spatial axis $x=0$, aligning with the expected characteristics of the reference solution. The reason may be that the map $F$ is learned from historical data and will preserve spatial symmetry. In contrast, the model state estimated by EnKF will be influenced by model error, resulting in the loss of this  symmetry.
	
	We describe the filtering distribution $p(u_k|Y_k)$ for the test data through KDE in Fig.\;$\ref{ac_2}$. The confidence intervals of $u_k$ estimated by the proposed method keep symmetric about the spatial axis $x=0$ due to the characteristics of $F$. In contrast, with EnKF, the first 20 components are estimated more accurately due to the availability of observation containing information pertaining to these components. The remaining components, which lack observational data, are estimated with higher uncertainty. This indicates that our proposed approach, employing a data-driven method, retains the physical characteristics of the state space even when observational information is limited.

	To investigate  the impact of different operators, we estimate the model state $u_k$ employing different operators for the reconstruction map: NNs, L-C operator and L-L operator in Fig.\;$\ref{ac_3}$. The resulting RMSEs and spreads are presented in Table $\ref{ac_4}$, which shows that NNs has the best accuracy and robustness. L-C operator, as a non-parametric approach, yields smoother estimations of the model state. In contrast, L-L operator exhibits strong sensitivity to small perturbations due to its use of weighted linear regression, resulting in more oscillatory estimations of the model state.
	
	\begin{table}[h]
		\centering
		\begin{tabular}{lllll}
			\hline
			& Neural Networks & L-C operator & L-L operator & EnKF     \\ \hline
			RMSE   & 0.146922        & 0.181342   & 0.202467   & 0.177576 \\
			Spread & 0.067872        & 0.067539   & 0.314727   & 0.178170 \\ \hline
		\end{tabular}
		\caption{RMSEs and spreads of test data by  different operators for the reconstruction map: NNs, L-C operator and L-L operator.}
		\label{ac_4}
	\end{table}
	
	\section{Conclusions}\label{section6}
	In this paper, we presented  DMD-T EnKF for online estimation without the underlying dynamics. In the offline stage, we constructed the transition model for observation  and the reconstruction map. The process of estimating the transition model can be interpreted  by  the perspective of EM algorithm. Then the reconstruction map was learned with the embedding and its corresponding model state. In the online stage, we utilized the surrogate model to estimate the  model state. Furthermore, we developed this framework with nonparametric time series prediction methods, namely KNN-T EnKF. The transition model and the reconstruction map has been estimated with analog forecasting methods, without learning the specific form. The  numerical results showed that the proposed approach  can use data pairs of the dynamical model and  accurately approximate the filtering distribution. The stability and accuracy can be enhanced by Takens theorem.

\section*{Aknowledge}
L. Jiang acknowledges the support of NSFC 12271408.

\smallskip
\bigskip
\textbf{Conflict of interest statement: }
The authors have no conflicts of interest to declare. All co-authors have seen and agree with the contents of the manuscript.

	\bibliographystyle{elsarticle-num}
	\bibliography{refs}
	\end{document}